\documentclass[11pt]{article}

\usepackage[colorlinks=true, allcolors=blue]{hyperref}

\usepackage{abstract, algorithm, algpseudocode, amsmath, amsthm, amssymb, appendix, bm, booktabs, caption, cite, etoolbox, lettrine, listofitems, listings, lmodern, fancyhdr, graphicx, hyperref, multirow, neuralnetwork, subfigure, tikz, titlesec, titling, wrapfig, xcolor}

\newcommand{\dropcap}[1]{\lettrine[lines=2,lraise=0.05,findent=0.1em, nindent=0em]{{\sffamily{#1}}}{}}

\usepackage[english]{babel}
\usetikzlibrary{shapes.geometric, arrows, bending}

\definecolor{gray}{cmyk}{0.05,0.1,0.1,0}

\algnewcommand{\Inputs}[1]{%
  \State \textbf{Inputs:}
  \Statex \hspace*{\algorithmicindent}\parbox[t]{.8\linewidth}{\raggedright #1}
}

\usepackage[margin=0.5in, top=1.2in, bottom=1in]{geometry}
\pretitle{
    \begin{center}
    \vspace{-15mm} 
    \Huge
    \bfseries 
    \sffamily
} 
\posttitle{
    \end{center}
    \vspace{1mm}
} 
\renewcommand\thesection{\Roman{section}} % Roman numerals for the sections
\renewcommand\thesubsection{\roman{subsection}} % roman numerals for subsections
\titleformat{\section}[block]{\Large\sffamily\bfseries}{\thesection.}{1em}{} % Change the look of the section titles
\titleformat{\subsection}[block]{\large\bfseries\sffamily}{\thesubsection.}{1em}{} % Change the look of the section titles

\title{Neural parameter calibration for large-scale multi-agent models} 

\date{}

\author{
	\bfseries{Thomas Gaskin$^{[1, \star]}$}  
\and
	\bfseries{Grigorios A. Pavliotis$^{[2]}$} 
\and 
	\bfseries{Mark Girolami$^{[3, 4]}$}
}

\begin{document}

\captionsetup{labelfont=bf, font=small}

\maketitle

\vspace{-10mm}
\begin{changemargin}{27pt}{27pt}
{\scriptsize\centering
$^\mathbf{1}$Department of Applied Mathematics and Theoretical Physics, University of Cambridge, Cambridge CB3 0WA, United Kingdom; $^\mathbf{2}$Department of Mathematics, Imperial College London, London SW7 2AZ, United Kingdom; $^\mathbf{3}$Department of Engineering, University of Cambridge, Cambridge CB2 1PZ, United Kingdom; $^\mathbf{4}$The Alan Turing Institute, London NW1 2DB, United Kingdom 

\medskip $^\star$To whom correspondence should be addressed: trg34@cam.ac.uk

}
\end{changemargin}

\vspace{6mm}
	
\begin{abstract} 

\noindent Computational models have become a powerful tool in the quantitative sciences to understand the behaviour of complex systems that evolve in time. However, they often contain a potentially large number of free parameters whose values cannot be obtained from theory but need to be inferred from data. This is especially the case for models in the social sciences, economics, or computational epidemiology. Yet many current parameter estimation methods are mathematically involved and computationally slow to run. In this paper we present a computationally simple and fast method to retrieve accurate probability densities for model parameters using neural differential equations. We present a pipeline comprising multi-agent models acting as forward solvers for systems of ordinary or stochastic differential equations, and a neural network to then extract parameters from the data generated by the model. The two combined create a powerful tool that can quickly estimate densities on model parameters, even for very large systems. We demonstrate the method on synthetic time series data of the SIR model of the spread of infection, and perform an in-depth analysis of the Harris-Wilson model of economic activity on a network, representing a non-convex problem. For the latter, we apply our method both to synthetic data and to data of economic activity across Greater London. We find that our method calibrates the model orders of magnitude more accurately than a previous study of the same dataset using classical techniques, while running between 195 and 390 times faster.

\bigskip

\end{abstract}
\vspace{6mm}
\hrule
\vspace{1mm}
\setcounter{tocdepth}{1}
\tableofcontents
\vspace{4mm}

\newpage

\twocolumn

\section{Introduction}
\vspace{2cm}

\dropcap{W}e live in an age of complexity. Thanks to an array of sophisticated and potent computational resources, paired with vast and fruitful reservoirs of data, researchers can increasingly see social, economic, biological, or epidemiological processes as the complex, self-organising, and dynamical systems they are. A line of great current interest in the mathematical sciences is to calibrate parameters of mathematical models using data, thereby creating computationally efficient, theoretically grounded, and socially beneficial predictive tools. For instance, models of the spread of contagion continues to inform much government policy during the ongoing COVID-19 pandemic \cite{Verity_2020, Hogan_2020, Flaxman_2020, Maier_2021, IHME_2020}: predictions are generated, subsequently compared to live data, and the underlying model revised accordingly. Computational models have also been used to understand the dynamics of crime and urban violence \cite{Bhavnani_2013, Helbing_2014}, pedestrian dynamics \cite{Helbing_1995}, synchronised oscillations \cite{Kuramoto_1975}, social network dynamics \cite{Castellano_2009, DelVicario2016, Bak_Coleman_2022}, chemotaxis and flocking \cite{Keller_1971, Vicsek_1995, CuckerSmale07}, population dynamics \cite{Arditi_1989, Winkelmann_2021}, systemic risk \cite{Giesecke2019}, and the dynamics of economic systems \cite{HarrisWilson78, Ellam_2018}. 

Two modelling paradigms have established themselves over the years: in the first, a system of coupled differential equations is iteratively solved to simulate the behaviour of a \emph{finite} number of interacting particles. This is at times implemented using what are called \emph{agent-based models} (ABMs): a collection of entities (\emph{agents}) moving through space and time, interacting with each other and their environment, and adapting their behaviour or even learning in the process. In the second approach, a model is realised as a discretised version of equations in \emph{continuous time and space} \cite{Martinez_1992, Ajelli_2010}; this approach opens models up to analysis using mean-field theory, stochastic analysis, statistical physics, and kinetic theory \cite{Lighthill_1955, Toscani_2006, Wang_2017, Carrillo_2019}. A large and sophisticated toolset of statistical methods exists to estimate parameters, interaction kernels, or network structures from data, such as maximum likelihood estimators \cite{Kasonga_1990, Pavliotis_2007, Chen_2021, Sharrock_2021, Liu_Qiao_2022}, Markov-Chain Monte Carlo methods based on a Bayesian paradigm \cite{Hastings_1970, Kaipio_2006, Stuart_2010, Gelman_2013}, martingale estimators \cite{Pavliotis_2021}, estimation of active terms in ODE and PDE systems (see \cite{Ramsay_2017} for a review), entropy maximisation \cite{Timme_2014}, and regression-based learning methods \cite{Lu_2021a, Lu_2021b}. More recently, a promising new method has emerged in the form of artificial neural nets. Neural networks have of course prominently been used as powerful pattern-recognition devices and predictive models \cite{Wei_2022}, but, as they become more and more accessible to the scientific community at large, researchers are beginning to apply their computational capabilities across the mathematical disciplines, including in fields heretofore dominated by more classical methods: examples include finding solutions of partial differential equations \cite{Raissi_2019, Kharazmi_2021} or parameter estimation of multi-agent models \cite{Goettlich_2021, Dyer_2022}. Neural networks, and especially deep neural networks, are mathematically little understood, and their theoretical underpinnings sparse and mainly restricted to shallow networks (networks with only one hidden layer) \cite{Sirignano_2018, Sirignano_2020}. This major drawback notwithstanding, recent results seem to indicate that their computational performance can often outstrip that of other, thus far more rigorously understood approaches, though the method still lies in its infancy.

This work is a contribution to the general push to better understand the possibilities of neural nets as calibration tools for mathematical models of complex dynamics. We present and investigate a simple yet powerful computational scheme to obtain probability densities for model parameters from data. The method combines classical numerical models with machine learning, and in the following case studies we recover probability densities from noiseless and noisy, synthetic and real, and time series and steady-state data. The case of a non-convex problem is also considered. Using the well-known SIR model of contagious diseases, we estimate parameter densities from a time series modelling the diffusion of infection through a population on a two-dimensional domain. In a second study, we use the Harris-Wilson model of economic activity on a network \cite{HarrisWilson78} to learn parameter densities both from synthetic steady-state data as well as a real dataset of activity across Greater London. In doing so we revisit an earlier study of the same dataset that used Bayesian methods to estimate parameters \cite{Ellam_2018}. Our proposed method can find estimates for model parameters in seconds, even for large systems, and provides parameter densities for the London dataset in nearly one minute where the Bayesian approach took between 3 and 7 hours. At the same time, the quality of the calibration (in terms of prediction error) is improved by two to three orders of magnitude.

The scope of our approach covers models formulated as coupled differential equations of the kind  
\begin{equation}
	\mathrm{d} \boldsymbol{\varphi} = f(\boldsymbol{\varphi}; \mathbf{x}, t, \boldsymbol{\lambda}) \mathrm{d}t,
	\label{eq:ODE_system_1}
\end{equation}
or, in the stochastic case, 
\begin{equation}
	\mathrm{d}\boldsymbol{\varphi} = f(\boldsymbol{\varphi}, \mathbf{x}, t, \boldsymbol{\lambda})\mathrm{d}t + g(\boldsymbol{\varphi}, t)\mathrm{d}\mathbf{B}_t,
	\label{eq:SDE_system_1}
\end{equation}
where $\boldsymbol{\varphi} \in \mathbb{R}^N$ is the state vector, $\mathbf{x}$ is a space-like variable, $\boldsymbol{\lambda} := (\lambda_1, ..., \lambda_p) \in \mathbb{R}^p$ a set of scalar parameters, and $\mathbf{B}$ an $N$-dimensional stochastic process (such as a Wiener process). (We include $\mathbf{x}$ to allow for infinite-dimensional problems leading to (stochastic) partial differential equations, though in practice discretisation will often lead to it being absorbed into the state vector.) In a \emph{neural} ODE or SDE, the scalar parameters $\boldsymbol{\lambda}$ are the outputs of a \emph{neural network}, whose \emph{internal parameters} are to be learned from data \cite{Kidger_2022}. In this work, we investigate the use of such neural ODEs as calibration tools. 

In many cases, the spatial topology is given by a network structure, such that the dependency on $\mathbf{x}$ is more specifically a dependency on a \emph{graph adjacency matrix} $\mathbf{A}$; this is found in many contemporary models of dynamical systems. One might also like to learn network structures from data. By vectorising $\mathbf{A}$, it can be conceived of as a single vector of parameters $\boldsymbol{\lambda}$, transforming the question into one that can indeed be considered within our proposed framework. However, as the scale of this problem is typically of a different order of magnitude, we address it in future work.

The dynamics of systems governed by equations such as eq. [\ref{eq:ODE_system_1}] can depend sensitively on the choice of parameters $\boldsymbol{\lambda}$. For mathematical models such as those used in the social sciences or computational biology, theoretical estimates for these parameters are often difficult obtain: what is the reproduction number of a novel disease, and how susceptible are different age groups to the disease? What is a good model for the topology of a social network? How should we estimate return on capital rates for different agents in economic models? What are recovery rates for different animal species subject to predator-prey dynamics? Such parameters are difficult to measure, and should instead be extracted from data.

\begin{figure}[!t]
\hrule{}
\begin{minipage}{.5\textwidth}
   \begin{tikzpicture}[->,scale=3.5, 
   node1/.style={minimum height=1cm, text centered, font=\bfseries, text width=3.4cm},
   node2/.style={minimum height=1cm, text centered, font=\bfseries, text width=2.1cm}]
   \node [node1] (NN) at (90:1cm)  {Neural Net};
   \node [node1] (Params) at (0:1cm) {Parameter estimates $\boldsymbol{\hat{\lambda}}$};
   \node [node1] (ABM) at (-90:1cm) {Numerical Solver (e.g. ABM)};
   \node [node2] (predictions) at (-170:1cm) {Predicted data $\hat{\mathbf{T}}$};
   \node [node2] (data) at (-190:1cm) {Observed data $\mathbf{T}$};
   \node [node1] (loss) at (0:0cm) {Loss \\ functional $J$};
   \node (jct) at (0:-0.5cm) {};
   \draw (70:1cm)  arc (70:10:1cm) node[midway, fill=white, sloped]{\normalsize{outputs}};
   \draw (-10:1cm) arc (-10:-60:1cm) node[midway, fill=white, sloped, rotate=0] {\normalsize{used to run}};
   \draw (-120:1cm) arc (-120:-160:1cm) node[midway, fill=white, sloped] {\normalsize{produces}}; 
   \draw (-200:1cm) arc (-200:-250:1cm) node[midway, fill=white, sloped] {\normalsize{input to}};
   \draw (loss.north) -- (NN.south) node[midway, fill=white] {\normalsize{used to train}};
   \draw[-] (predictions.east) --  (jct.west);
   \draw[-] (data.east) -- (jct.west);
   \draw (jct.west) -- (0:-0.4);
   \end{tikzpicture}
\end{minipage}   
\begin{minipage}{.5\textwidth}

\begin{algorithm}[H]
\caption{Single training epoch}
\begin{algorithmic}[1]
    \Inputs{$\mathbf{T} = (\boldsymbol{\varphi}_1, ... , \boldsymbol{\varphi}_L)$ (time series of length $L$) \\ $B \leq L$ (batch size)}
    \For {$t \in \{1, ..., L-B\}$}
    \State {$\boldsymbol{\hat{\lambda}} \gets u_\theta(\boldsymbol{\varphi}_t, ..., \boldsymbol{\varphi}_{t+q})$ (parameter estimates)}
    \State {$\hat{\boldsymbol{\varphi}}_0 \gets \boldsymbol{\varphi}_t$}
    \For {$b \in \{1, ..., B\}$}
    \State $\hat{\boldsymbol{\varphi}}_b \gets \int f(\hat{\boldsymbol{\varphi}}_{b-1}; \boldsymbol{\hat{\lambda}})\mathrm{d}t$
    \EndFor
    \State Calculate loss $J(\hat{\boldsymbol{\varphi}}_1, ..., \hat{\boldsymbol{\varphi}}_{B}, \boldsymbol{\varphi}_{t+1}, ..., \boldsymbol{\varphi}_{t+B})$
    \State Calculate gradient $\nabla_{\boldsymbol{\theta}} J$
    \State Update $\boldsymbol{\theta}$ using backpropagation
    \EndFor
\end{algorithmic}
\end{algorithm}

\end{minipage}
\caption{The methodological pipeline proposed in this work. The neural net $u_\theta$ takes $q$ time series elements as input and outputs parameter predictions. These predictions are fed into a numerical solver, which produces a predicted time series. The true and predicted time series are used to generate a loss functional, which in turn can be used to train the neural net's internal parameters $\boldsymbol{\theta}$. The goal is to retrieve the true parameters $\boldsymbol{\lambda}$ from the data. A single pass over the entire dataset is called an \emph{epoch}. The dataset is processed in \emph{batches}, meaning the loss is calculated over $B$ steps of the time series before the weights are updated. If $B=L$, training is equivalent to batch gradient descent; if $B=1$, training is equivalent to stochastic gradient descent. The integral in line 6 represents an arbitrary numerical scheme to solve eq. [\ref{eq:ODE_system_1}].}
\label{diag:method}
\end{figure}

\section{Methodology}
\label{section:methodology}
\subsection{Obtaining probability densities from Neural equations} We present a method to estimate parameter densities of ODE or SDE systems by training a neural net to find a set of parameters $\boldsymbol{\hat{\lambda}}$ that, when inserted into the model equations eq. [\ref{eq:ODE_system_1}], reproduce a given time series $\mathbf{T} = (\boldsymbol{\varphi}_1, ..., \boldsymbol{\varphi}_L)$. A \emph{neural network} is a function $u_\theta: \mathbb{R}^{N \times q} \to \mathbb{R}^p$, where $q \geq 1$ represents the number of time series steps that are passed as input (cf. SI appendix fig. \ref{diag:neural_net}). Its output are the estimated parameters $\boldsymbol{\hat{\lambda}}$, which are used to run a numerical solver for $B$ iterations ($B$ is the \emph{batch size}) to produce an \emph{estimated time series} $\mathbf{\hat{T}}(\boldsymbol{\hat{\lambda}}) = (\hat{\boldsymbol{\varphi}}_i, ..., \hat{\boldsymbol{\varphi}}_{i+B})$. This in turn is used to train the internal parameters $\boldsymbol{\theta}$ of the neural net (the \emph{weights and biases}) via a \emph{loss functional} $J(\hat{\mathbf{T}}, \mathbf{T})$. A common choice for $J$ is the $l^2$ norm, which we will use in this work. As $\boldsymbol{\hat{\lambda}} = \boldsymbol{\hat{\lambda}}(\boldsymbol{\theta})$, we may calculate the gradient $\nabla_{\boldsymbol{\theta}}J$ and use it to optimise the internal parameters of the neural net using a backpropagation method of choice; popular choices include stochastic gradient descent, Nesterov schemes, or the Adam optimizer \cite{Nakkiran_2021, Kingma_2014}. Calculating $\nabla_{\boldsymbol{\theta}}J$ thus requires differentiating the predicted time series $\mathbf{\hat{T}}$, and thereby the system equations [\ref{eq:ODE_system_1}], with respect to $\boldsymbol{\hat{\lambda}}$. In other words: the loss function contains knowledge of the dynamics of the model. Finally, the true data is once again input to the neural net to produce a new parameter estimate $\hat{\boldsymbol{\lambda}}$, and the cycle starts afresh. A single pass over the entire dataset is called an \emph{epoch} (cf. fig. \ref{diag:method}). In the following we will always let $q=1$. The technicalities of the differentiation procedure are handled by the auto-differentiation scheme, which treats random vectors as constants. For non-differentiable convex functions (such as $\Vert \cdot \Vert$), the subgradient of minimum norm is used, see e.g. \cite{pytorch_autograd}.

This solves an optimisation problem, but it does not provide us with any sort of confidence intervals for the predictions. To obtain probability densities, we exploit the fact that as the model trains, it traverses the parameter space $\mathbb{R}^p$, calculating a loss value $J$ at each estimate $\boldsymbol{\hat{\lambda}}$. This produces a \emph{loss potential} $J(\boldsymbol{\hat{\lambda}}): \mathbb{R}^p \to \mathbb{R}$, from which a posterior density can be estimated via 
\begin{equation}
    \pi(\hat{\boldsymbol{\lambda}} \vert \mathbf{T}) \sim \exp(-J(\hat{\mathbf{T}}, \mathbf{T})) \pi^0(\hat{\boldsymbol{\lambda}}),
\end{equation}
with $\pi^0$ the prior. The marginals are then proportional to
\begin{equation}
	\rho(\lambda_i) \sim \int \exp(-J) \mathrm{d}\hat{\boldsymbol{\lambda}}_{-i}
	\label{eq:marginal_density}
\end{equation}
with the subscript $-i$ signifying that we are not integrating over $\hat{\lambda}_i$. In the following, we initialise the neural net's internal parameters such that the prior is a uniform density. Non-convexity can be dealt with by training the neural net multiple times on the same dataset, each time with a different \emph{initialisation}, thereby increasing the volume of the sampled parameter space. This comes at a computational cost, but multiple runs can be parallelised and run concurrently, thereby greatly improving performance.

Both the noiseless (eq. [\ref{eq:ODE_system_1}]) and the noisy (eq. [\ref{eq:SDE_system_1}]) version of the equations can be used when running the numerical solver (i.e. the forward pass of the pipeline alg. \ref{diag:method}). Running the solver without noise circumvents having to make any assumptions about the nature of the true underlying noise in the dataset, which will often be unknown. The neural net will then simply fit the best noiseless model to the dataset. However, noise can be added to the solver when randomness is an inherent part of the dynamics, and can help make the predictions more robust. Alternatively, the variance of the noise can itself be learned. All these scenarios will be demonstrated in this work.

\subsection{Using Neural equations in practice} It is important to reiterate that the neural net is \emph{not}  fitting a dataset in the traditional sense, but rather producing a set of parameters that, when plugged into the governing equations, \emph{generate} the dataset (or an estimate thereof). Differentiating the loss function (and thereby the physical equations) may seem daunting; but many of the standard machine learning libraries have auto-differentiation features, and so the differentiation procedure need not (and if possible should not) be implemented manually. In practice, the bottleneck of our method will lie in writing a fast numerical solver that is also compatible with the differentiation procedure of the machine learning package used. Explicitly iterating over agents or network nodes should be avoided, and the dynamics instead be implemented as operations on the entire state vector $\boldsymbol{\varphi}$. Pre-implemented matrix and vector operations should be used as much as possible, as these typically (1) use efficient, tested, and pre-compiled algorithms, and (2) are compatible with auto-differentiation features. We have implemented an open-source code package (see below), written such that it can be easily extended and adapted to further models. See the supplementary material and the README files in the repository.

Aside from computational considerations, when learning parameters, several fundamental mathematical questions must be carefully considered. \textbf{(1)} First, \emph{the model's dynamical range should be analysed} in order to ascertain whether there are non-identifiable regimes. Dynamical systems may gravitate toward attractors and steady equilibria which are independent of certain model parameters and thus cannot be used to train the net. \textbf{(2)} \emph{Numerical stability} must be maintained throughout the training process, e.g. by rescaling the neural net output to prevent numerical overflow. In certain parameter regimes, dynamical systems may also exhibit chaotic behaviour that inhibit the learning process. The neural net must be kept from straying into the `danger zones' of potential numerical instability as it traverses the parameter space during training. \textbf{(3)} To speed up computation, a suitable \emph{neural net architecture} should be chosen that encapsulates as much information on the parameters as can reasonably be assumed a priori. For instance, using activation functions that guarantee parameters remain in valid ranges may be conducive: if parameters are probabilities, an activation function that maps into $[0,1]$ (e.g. a sigmoid) may be an appropriate choice for the final layer. 
\subsection{Code and data availability} All code and data can be found under \url{https://github.com/ThGaskin/NeuralABM}. It is easily adaptable to new models and ideas. The code uses the \texttt{utopya} package\footnote{\href{https://utopia-project.org}{utopia-project.org}, \href{https://utopya.readthedocs.io/en/latest/}{utopya.readthedocs.io/en/latest}} \cite{Riedel2020, Sevinchan_2020} to handle simulation configuration and efficiently read, write, analyse, and evaluate data. This means that the model can be run by modifying simple and intuitive configuration files, without touching code. Multiple training runs and parameter sweeps are automatically parallelised. The neural core is implemented using \texttt{pytorch}\footnote{\href{https://pytorch.org}{pytorch.org}}. All datasets used in this work, including the synthetic data, have been made available, together with the configuration files needed to reproduce the plots. Detailed instructions are provided in the supplementary material and the repository.

\section{Application to time series data: diffusive SIR model of epidemics}
\label{section:SIR}
\begin{figure*}
    \begin{minipage}{1\textwidth}		
    \includegraphics[width=0.245\textwidth]{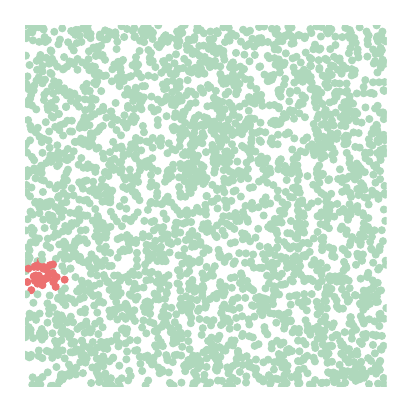}
    \includegraphics[width=0.245\textwidth]{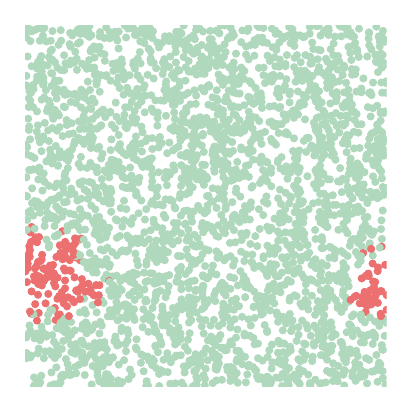}
    \includegraphics[width=0.245\textwidth]{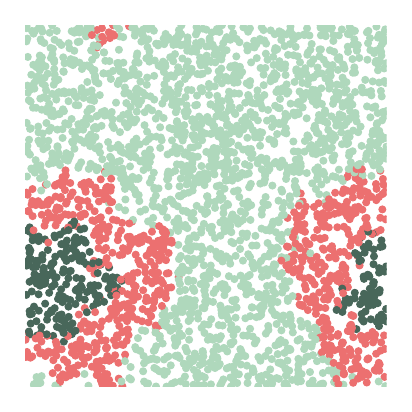}
    \includegraphics[width=0.245\textwidth]{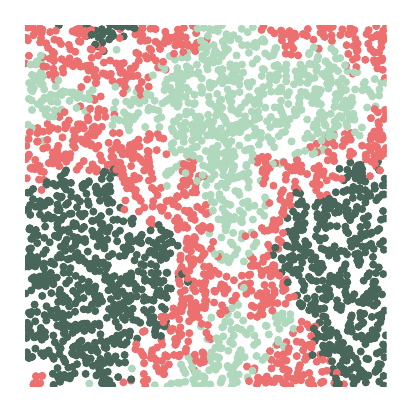}
\end{minipage}		
\caption{Diffusion of infection (red) through the agent population on a two-dimensional domain with periodic boundary conditions. Dark green: recovered agents.}
\label{fig:SIR_frames}
\end{figure*}
We first demonstrate our methodology by applying it to time series data generated by a classical agent-based model of epidemics. Consider $N$ agents moving around a square domain $[0, L]^2$ with periodic boundary conditions, $0 < L \in \mathbb{R}$; each agent has a \emph{position} $\mathbf{x}_i$, and a \emph{state} $k_i \in \{ \text{S, I, R} \}$. All agents with $k_i = \text{S}$ are \emph{susceptible} to the disease. If a susceptible agent lies within the \emph{infection radius} $r$ of an \emph{infected agent} (an agent with $k_i = \text{I}$), they are infected with \emph{infection probability} $p$. After a certain \emph{recovery time} $\tau$, agents \emph{recover} from the disease (upon which $k_i = \text{R}$); each agent's time since infection is stored in a state $\tau_i$. Agents move randomly around the space with diffusivities $\sigma_S, \sigma_I$, and $\sigma_R$. Each iteration of the agent-based model thus consists of the following steps:
\begin{algorithm}[H]
\caption{Single iteration of the SIR model}\label{alg:SIR_ABM}
\begin{algorithmic}[1]
\For {all agents $j$ with $k_j = \text{I}$}
\For {all agents $i$ with $k_i = \text{S}$}
\If {$d(\mathbf{x}_i, \mathbf{x}_j) \leq r$}
   \State $k_i \gets \text{I with probability } p$
   \State $\tau_i \gets 0$
\EndIf
\EndFor
\State $\tau_j \gets \tau_j +1$
\EndFor
\For {all agents $i$ with $k_i = \text{I}$}
	\If {$\tau_i \geq \tau$}
		\State $k_i \gets \text{R}$
	\EndIf
\EndFor
\For {all agents}
	\State Move agent randomly with respective diffusivity
\EndFor
\end{algorithmic}
\end{algorithm}
\noindent The function $d$ is the distance metric on the torus
\begin{equation*}
    d(\mathbf{x}, \mathbf{y})^2 = \sum_i \min \left(\vert x_i - y_i\vert, L-\vert x_i - y_i \vert \right)^2,
\end{equation*}
and we initialise the ABM with a single infected agent at a random location. We set the infection time to $\tau=14$, the infection probability to $p=0.2$, and the infection radius to $r=0.3$. The space has dimension $L=10$, and the generated time series contains 100 time steps. 

\begin{figure}[!hb]
    \includegraphics[width=8.7cm]{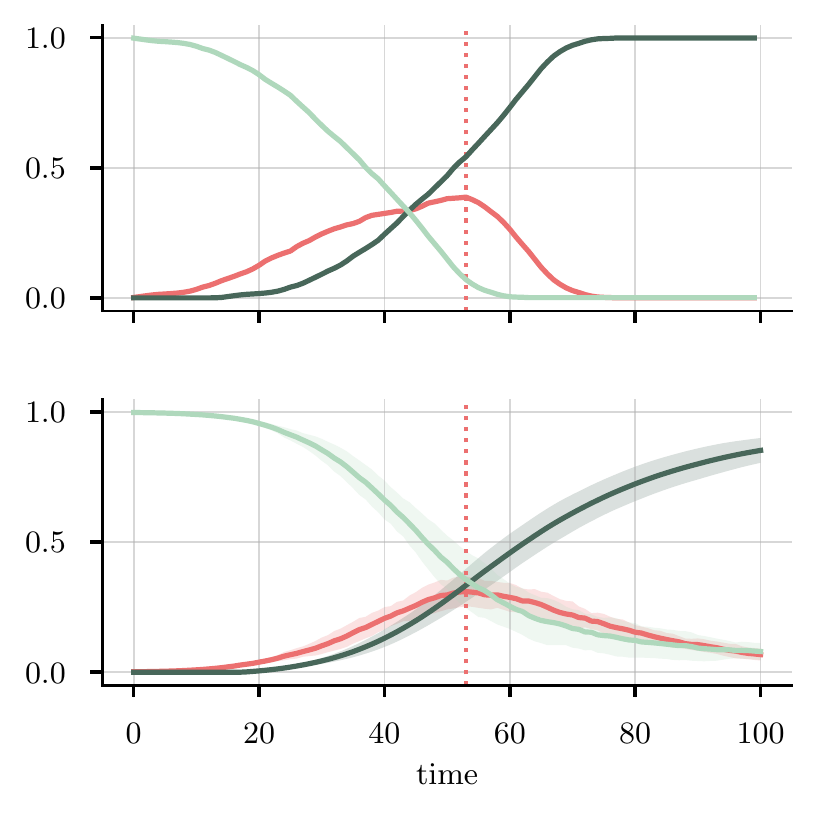}
    \caption{SIR densities $\text{S}(t)$ (light green), $\text{I}(t)$ (red), and $\text{R}(t)$. Top: sample densities used to train the neural net, generated by running the agent-based model with $N=3000$ agents for 100 iterations. Bottom: predicted densities using the neural net outputs after training, averaged over 20 different initialisations, with the errorbands showing one standard deviation. For each initialisation, the neural net is trained for 70 epochs, with a batch size of 90. Dotted lines: times of peak infection. Total runtime: 44.8s.}
    \label{fig:SIR_densities}
    \vspace{2cm}
\end{figure}
\begin{figure*}[ht!]
\centering
    \includegraphics[width=0.32\textwidth]{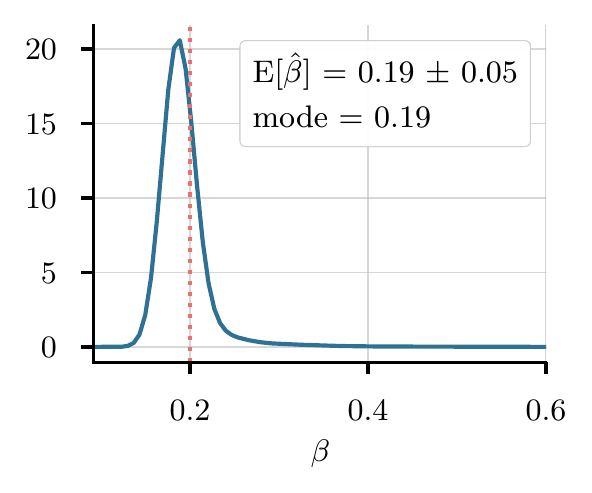}
    \includegraphics[width=0.32\textwidth]{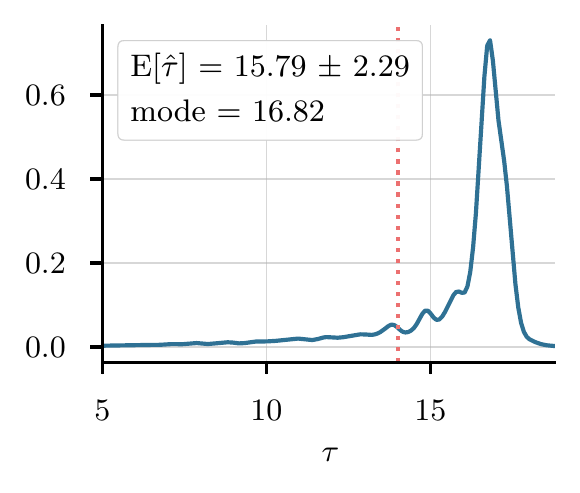}
    \includegraphics[width=0.32\textwidth]{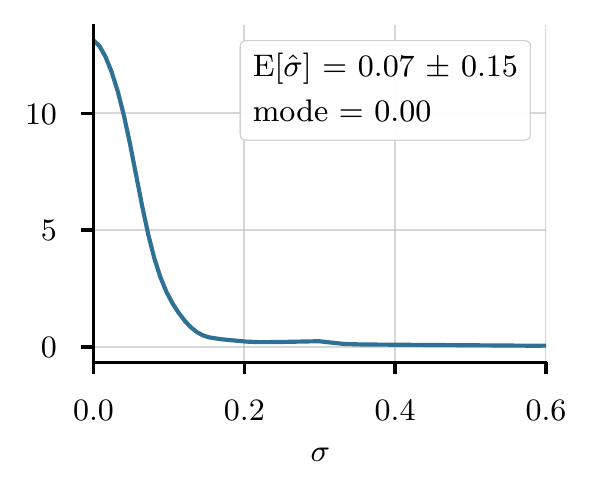}
\caption{Marginal densities for the parameters $\boldsymbol{\lambda} = (\beta, \tau, \sigma)$, calculated using eq. [\ref{eq:marginal_density}], smoothed with a Gaussian kernel. The parameters used to run the agent-based model are $p=0.2$ and $\tau=14$ (indicated by the red dotted lines).}
\label{fig:SIR_predictions}
\end{figure*}

Let $S(\mathbf{x}, t)$ be the spatio-temporal distribution of susceptible agents (analogously $I$ and $R$). Assume we only observe the temporal densities 
\begin{equation*}
    \text{S}(t) = \dfrac{1}{N}\int_\Omega S(\mathbf{x}, t) \mathrm{d}\mathbf{x},
\end{equation*}
applicable to the spread of an epidemic where we only see the counts of infected and recovered patients without any location tracking or contact tracing. To these observations we now wish to fit the stochastic equations
\begin{gather*}
	\mathrm{d}\text{S} = - \beta \text{SI} \mathrm{d}t - \sigma \text{I} \circ \mathrm{d}W \nonumber \\
	\mathrm{d}\text{I} = (\beta \text{S} - \tau^{-1})\text{I}\mathrm{d}t +  \sigma \text{I} \circ \mathrm{d}W \nonumber \\
	\mathrm{d}\text{R} = \tau^{-1} \text{I}\mathrm{d}t,
\end{gather*}
where $W$ is a Wiener process, and $\circ$ represents the Stratonovich integral. We will do so by recovering the parameters $\boldsymbol{\lambda} = (\beta, \tau, \sigma) \in \mathbb{R}^3_+$ from observations of S, I, and R generated by the agent-based model. We expect $\beta \approx p$, as the likelihood of a susceptible agent coming in to contact with an infected agent approaches 1. Naturally, there is some data-model mismatch, which is intentional and in this case accounted for by the noise $\sigma$; in reality, $\sigma$ could for instance model errors in the estimates of the number of infected agents.

We use a shallow neural net with 20 neurons in the hidden layer and, since parameters are known to be positive and negative outputs would produce unpredictable behaviour, the modulus $x \mapsto \vert x \vert$ as an activation function. The loss function is the batch-averaged mean squared error,
\begin{equation}
J(\hat{\boldsymbol{\varphi}}_i, ..., \hat{\boldsymbol{\varphi}}_{i+B}, \boldsymbol{\varphi}_i, ..., \boldsymbol{\varphi}_{i+B}) = \dfrac{1}{B}\sum_{j=i}^{i+B} \Vert \hat{\varphi}_j - \varphi_j \Vert^2,
\label{eq:loss_function}
\end{equation}
and we use the Adam optimizer \cite{Kingma_2014} for the backpropagation. For a single initialisation we train the neural net for 70 epochs with a batch size of 90, and we run the model from 20 different initialisations.

Figure \ref{fig:SIR_predictions} shows the parameter predictions. We recover the infection probability with a maximum likelihood estimate of $\hat{\beta} = 0.19$, and a slightly overestimated infection time $\hat{\tau} = 16.82$. The most likely noise level is predicted to be 0, with an expectation value of $0.07 \pm 0.15$. In fig. \ref{fig:SIR_densities}, we use the neural net predictions from all 20 initialisations to calibrate the model. We see how, on average, the predicted time of peak infection matches the true time (dotted lines), though the density of infected agents both increases and decreases significantly more slowly than in the true dataset. Despite the significant data-model mismatch, the neural network therefore manages to make reasonable predictions — both in terms of the estimated parameters and the output of the calibrated ABM — in only a few seconds: training the model for all 20 initialisations took 45s on a standard laptop CPU.

\section{Application to a non-convex problem: the Harris-Wilson model of economic activity}
\label{section:Harris_Wilson}
In this section we analyse the connection between prediction uncertainty and noise in the training data, as well as comparing the method to more classical methods with regard to its predictive ability and computational performance. We do so using the Harris-Wilson model of economic activity on a network \cite{HarrisWilson78}, a non-convex problem for which the loss function has at least two global minima. In a first step, we will consider synthetic data, thereby avoiding any data-model mismatch, and giving us control over the variance in the data; thereafter we shall analyse a real dataset of economic activity in Greater London. Before presenting our results we briefly describe the model dynamics.

In the Harris-Wilson model, $N$ \emph{origin zones} are connected to $M$ \emph{destination zones} through a weighted, directed, \emph{complete bipartite network}, i.e. each origin zone is connected to every destination zone. Economic demand flows from the origin zones to the destination zones, which supply the demand. Such a model is applicable for instance to an urban setting, the origin zones representing e.g. residential areas, and the destination zones representing retail areas, shopping centres, or other areas of consumer activity. Let $\mathbf{C} \in \mathbb{R}^{N \times M}$ be the non-zero section of the full network adjacency matrix. The \emph{network weights} $c_{ij}$ quantify the convenience of travelling from origin zone $i$ to destination zone $j$: a low weight thus models a highly inconvenient route (e.g. due to a lack of public transport). Each origin zone has a fixed demand $O_i$. The resulting cumulative demand at some destination zone $j$ is given by 
\begin{equation*}
	D_j = \sum_{i=1}^{N} T_{ij},
\end{equation*}
$T_{ij}$ representing the \emph{flow of demand} from $i$ to $j$. The model assumption is that this flow depends both on the \emph{size} $W_j$ of the destination zone and the convenience of `getting~from~$i$~to~$j$':
\begin{equation*}
	T_{ij} = \dfrac{W_j^\alpha c_{ij}^\beta}{\sum_{k=1}^M W_k^\alpha c_{ik}^\beta} O_i.
\end{equation*}
The parameters $\alpha$ and $\beta$ represent the relative importance of size and convenience to the flow of demand from $i$ to $j$: high $\alpha$ means consumers value large destination zones (e.g. prefer larger shopping centres to smaller ones), high $\beta$ means consumers place a strong emphasis on convenient travel to destination zones. Finally, the sizes $W_j$ are governed by a system of $M$ coupled logistic equations:
\begin{equation}
	\mathrm{d}W_j = \epsilon W_j(D_j - \kappa W_j)\mathrm{d}t + \sigma W_j \circ \mathrm{d}B_j,
	\label{eq:HW_model}
\end{equation}
with given initial conditions $W_j(t=0) = W_{j, 0}$. Here, $\epsilon$ is a responsiveness parameter, representing the rate at which destination zones can adapt to fluctuations in demand, and $\kappa$ models the cost of maintaining a larger site per unit floor space (e.g. rent, utilities, etc.). We recognise the logistic nature of the equations: the change in size is proportional to the size itself, as well as to $D_j - \kappa W_j$. A low value of $\kappa$ favours larger destination zones (e.g. larger malls), a high cost favours smaller zones (e.g. local stores). 
In addition, the model eq. [\ref{eq:HW_model}] includes multiplicative noise with variance $\sigma \geq 0$, with $\circ$ signifying Stratonovich integration. This represents a perturbation of the net capacity term $D_j - \kappa W_j$ by a centred Gaussian. We thus have a system of $M$ coupled stochastic differential equations. In the noiseless case, the stable equilibrium is determined by
\begin{equation}
	\mathbf{W} = \kappa^{-1} \mathbf{D},
	\label{eq:steady_state}
\end{equation}
where $\mathbf{W} \in \mathbb{R}^M$ and $\mathbf{D} \in \mathbb{R}^M$ are the origin zone sizes and demands, respectively. The steady state is thus independent of the responsiveness $\epsilon$, which only affects the convergence rate to the equilibrium. We can therefore set $\epsilon=1$. Note also that $\alpha$ and $\beta$ are unitless, and thus unaffected by any scaling of the origin or destination zone sizes. $\kappa$ is given in units of cost/area, and scales accordingly. 

\begin{figure}[!b]
	\flushleft
	\begin{minipage}{\textwidth}
		\includegraphics[width=0.5\textwidth]{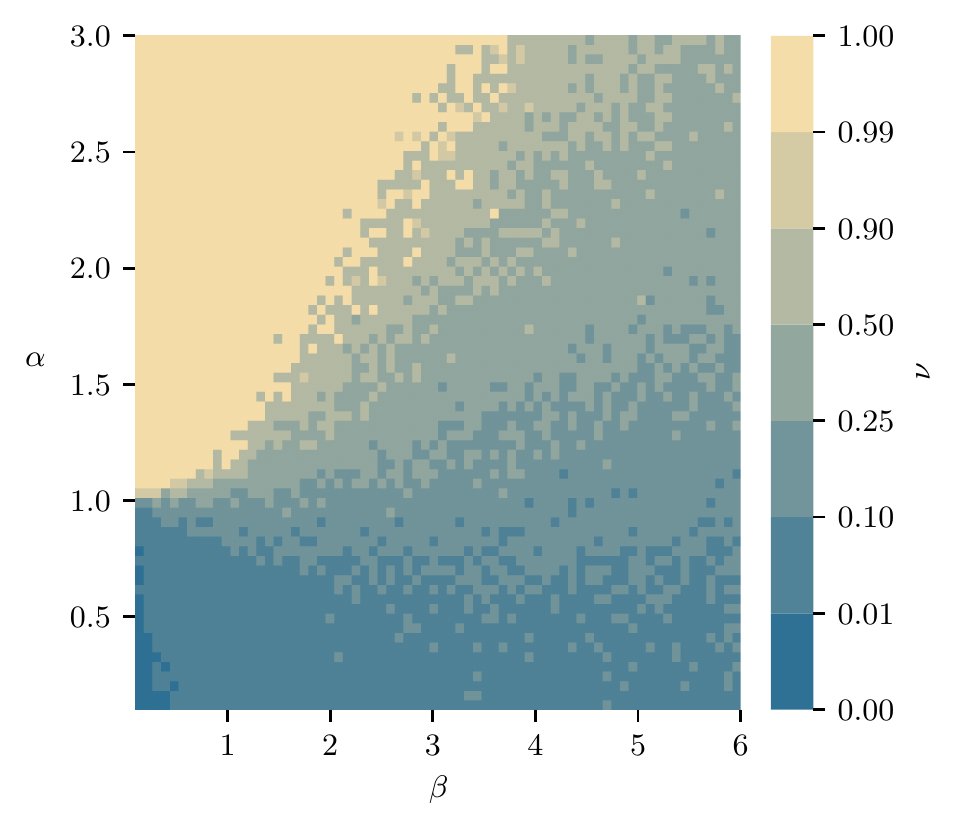}
	\end{minipage}
	\caption{The inequality parameter $\nu$ (eq. [\ref{eq:nu}]) of the destination zone sizes as a function of $\alpha$ and $\beta$ (holding $\kappa$ fixed). For high $\beta$ and low $\alpha$, a larger number of smaller centres can emerge, whereas for high $\alpha$ and relatively low $\beta$, a small number of super-centres can form, until finally the market becomes entirely dominated by a single zone.}
	\label{fig:inequality}
\end{figure}
In the stochastic case, the dynamics reach a steady-state equilibrium that is independent of the initial condition $\mathbf{W}_0$. A good indicator to assess the effect of the parameters $\alpha$ and $\beta$ on the system steady state is the \emph{inequality}
\begin{equation}
	\nu := \sup_{i, j}  \dfrac{W_j - W_i}{\sum_k W_k}  \in [0, 1];
	\label{eq:nu}
\end{equation}
$\nu=1$ indicates a completely monopolised market, while $\nu = 0$ indicates perfect equality of size, i.e. all zones are of equal size. Low values of $\alpha$ and high values of $\beta$ (i.e. low relative importance of zone size, high relative importance of convenience) lead to low overall inequality, due to a large collection of smaller zones emerging, all of roughly equal size (cf. fig. \ref{fig:inequality}). Conversely, low relative importance of travel convenience, and high relative importance of store size lead to the emergence of a small number of very large superstores, with most smaller centres dying out: this is reflected in the high values of $\nu$.

The steady state condition eq. [\ref{eq:steady_state}] is given by the $M$ coupled equations
\begin{equation}
	\sum_i \dfrac{c_{ij}^\beta O_i}{\sum_k W_k^\alpha c_{ik}^\beta} = \kappa W_j^{1-\alpha}, \ j = 1 ... M.
	\label{eq:steady_state_2}
\end{equation}
In the case of $\nu = 1$, eq. [\ref{eq:steady_state_2}] will be solved by any $\hat{\alpha}$ and $\hat{\beta}$, with $\hat{\kappa}$ uniquely given by
\begin{equation*}
	\hat{\kappa} = \dfrac{\sum_i O_i}{W_{\hat{k}}},
\end{equation*}  
where $\hat{k}$ is the index of the monopolising zone. This region is thus unlearnable in $\alpha$ and $\beta$ (with consumers only having a single option, both their preference for size and convenience is irrelevant). Similarly, the case $\nu = 0$ admits solutions that are independent of $\alpha$ (since all zones are of the same size, the size preference parameter becomes irrelevant). \emph{Parameters thus cannot be learned for} $\nu \in \{0, 1\}$, and throughout this work, we only consider datasets with $0 < \nu < 1$.

For any $\nu$, the triple $(\alpha=1, \beta=0, \kappa=\sum_i O_i / \sum_k W_k )$ represents a global minimum of the loss function, and if $0 < \nu < 1$, it can be shown that there exist \emph{at most} two global minima of the loss function (see \cite{HarrisWilson78} for a proof):

\medskip
\noindent \textbf{Fact.} \textit{In their stable equilibrium, the noiseless Harris-Wilson equations eq. [\ref{eq:steady_state}] for $0 < \nu < 1$ admit at most two solutions in $(\alpha, \beta, \kappa)$, one being the \emph{trivial solution} $(1, 0, \sum_i O_i / \sum_k W_k)$. For any $\nu$, the solution is unique in $\kappa$, where $\kappa$ is given by $\sum_i O_i/\sum_k W_k$.}

\subsection{Results on synthetic data}
We generate synthetic data of the steady-state sizes \hspace{\fill}\linebreak $\mathbf{W}^\star~:=~\mathbf{W}(t~\to~\infty)$, determined via $\sup_j \mathrm{d}W_j~<~\mathrm{tol}$ for some sufficiently small tolerance. We generate training data with different noise levels; the data has length $L=1$ in the noiseless and $L=4$ in the noisy case. We run the numerical solver without noise, since we do not wish to include any assumptions about $\sigma$ in the training process.

\begin{figure}[t!]
    \vspace{-5mm}
    \includegraphics[width=0.5\textwidth]{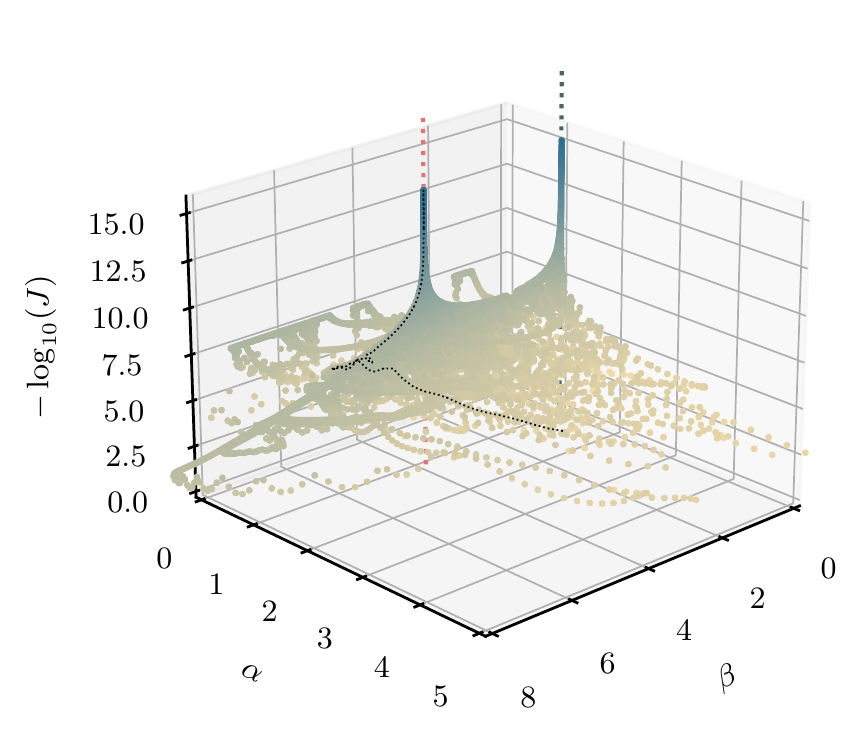}
	\caption{The negative logarithm of the loss potential $J$ as a function of the estimates $\hat{\alpha}$ and $\hat{\beta}$. Each point represents a single training step, and the model was trained for 10000 iterations over 100 different initialisations. The two minima of the potential at the true parameters ($\alpha=1.2, \beta =4, \kappa=2$, red line) and the trivial solution ($\alpha=1, \beta=0, \kappa=2$,  green line) are clearly visible. Also shown is an example trajectory the model takes as it converges to the (non-trivial) solution (black dotted line).}
	\label{fig:loss_landscape}
\end{figure}

Considering first the noiseless case, we wish to estimate a probability distribution for the model parameters $\{\alpha, \beta, \kappa \} \in \mathbb{R}_+^3$ with confidence bounds ($\epsilon$ is not learnable, as it does not affect the steady state). The steady state must meet the criterion $0 < \nu < 1$. We train the neural net 100 times from different initialisations for 10000 epochs each, performing a gradient descent step after each prediction (thus the batch size is $B=1$), and using the same network architecture as before. The loss function is the same as in the SIR case (eq. [\ref{eq:loss_function}]). Figure \ref{fig:loss_landscape} visualises the resulting loss potential in $\alpha$ and $\beta$. We clearly see the two global minima of $J$, one at the true parameters and one at the trivial minimum of the decoupled system $\alpha=1$, $\beta=0$. In running the model several times with different initialisations, we are thus able to deal with the non-convexity of the problem. 

\begin{figure}[hb!]
    \includegraphics[width=0.5\textwidth]{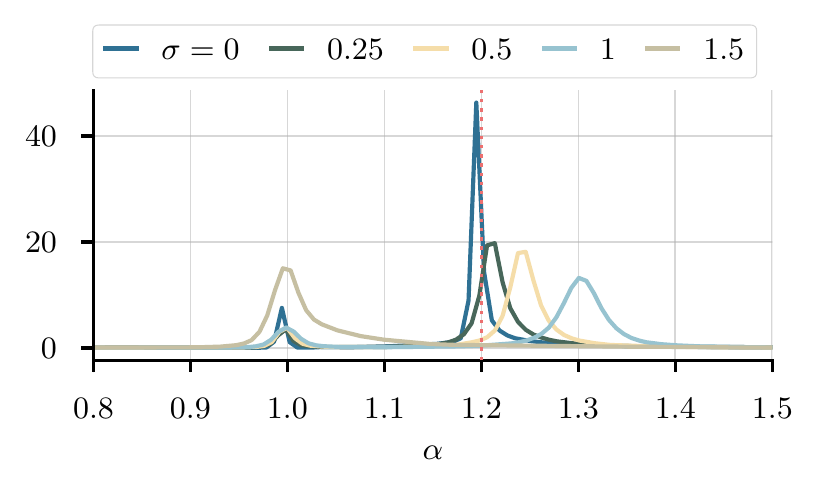}
    \caption{The marginal densities for $\alpha$ as a function of the noise in the data, smoothed with a Gaussian kernel. As we increase the noise, the peak width increases (see also fig. \ref{fig:HW_peak_widths}). Red dotted line: true value.  Similar results hold for the other parameters, see fig. \ref{fig:HW_additional_marginals} in the appendix.}
    \vspace{5mm}
    \label{fig:HW_synthetic_marginals}
	\includegraphics[width=0.5\textwidth]{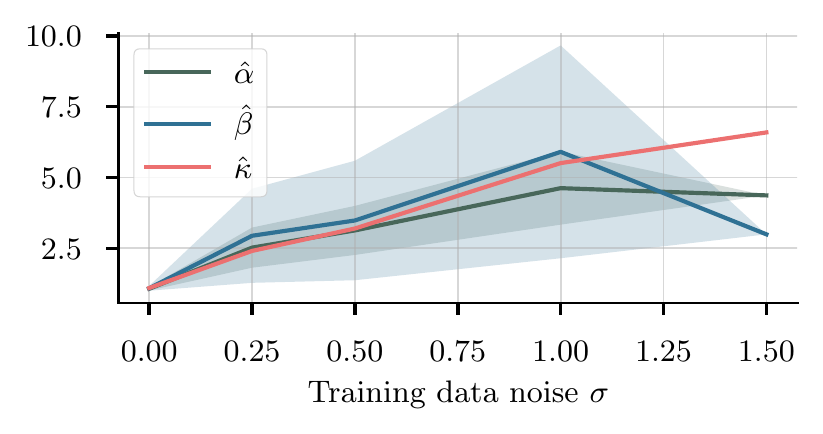}
	    \caption{Average peak width (defined as the width at half height) for the three parameters. Shaded area: standard deviation. For very high noise, only the trivial peak is found, wherefore the standard deviation becomes 0.}
    \label{fig:HW_peak_widths}
\end{figure}

Turning to the case of noisy data, fig. \ref{fig:HW_synthetic_marginals} shows the marginals for $\alpha$ as we increase the level of noise in the data. The noise is generated by a mean-zero, centred Gaussian, and its variance ranges from $0$ to a very high level of $\sigma=1.5$. In the noiseless case (dark blue) the non-trivial solution dominates, with a small peak at the trivial case $\alpha=1$. As the noise increases, so do the peak widths, until at $\sigma=1.5$ the non-trivial solution is no longer identiable, and the model collapses to the trivial value $\alpha=1$. Figure \ref{fig:HW_peak_widths} shows the average and standard deviation of the peak widths for all three parameters. Since $\kappa$ is uniquely identifiable, there is only one peak, hence the standard deviation is 0. For the other two, the disappearance of the non-trivial peak is evident at $\sigma=1.5$, where the standard deviation of the peak width becomes 0. For all three parameters, we observe an increase in the average peak width as the noise in the data increases. At $\sigma=0$, there is still some uncertainty in the predictions, originating from the uncertainty in the internal model weights $\boldsymbol{\theta}$.

Generating the 1 million samples in 3-dimensional space shown in figure \ref{fig:loss_landscape} takes about 15 minutes on a standard CPU. Far fewer data points would already have been sufficient to get good estimates of the parameters, since the model spends most of its time near the minima of the loss function. The time to run a single iteration increases with $(N+M)^2$, as the adjacency matrix of the network (which is not sparse) grows accordingly; the loss after a fixed number of iterations remains fairly constant (see SI appendix figure \ref{fig:HW_performance_analysis}). In the following section, we give a more detailed analysis of the model performance by comparing it to that of a classical Bayesian estimator.

\subsection{Results on data of economic activity in London}
\begin{figure*}[t!]
	\begin{minipage}[!t]{.64\textwidth}
		\includegraphics[width=\textwidth]{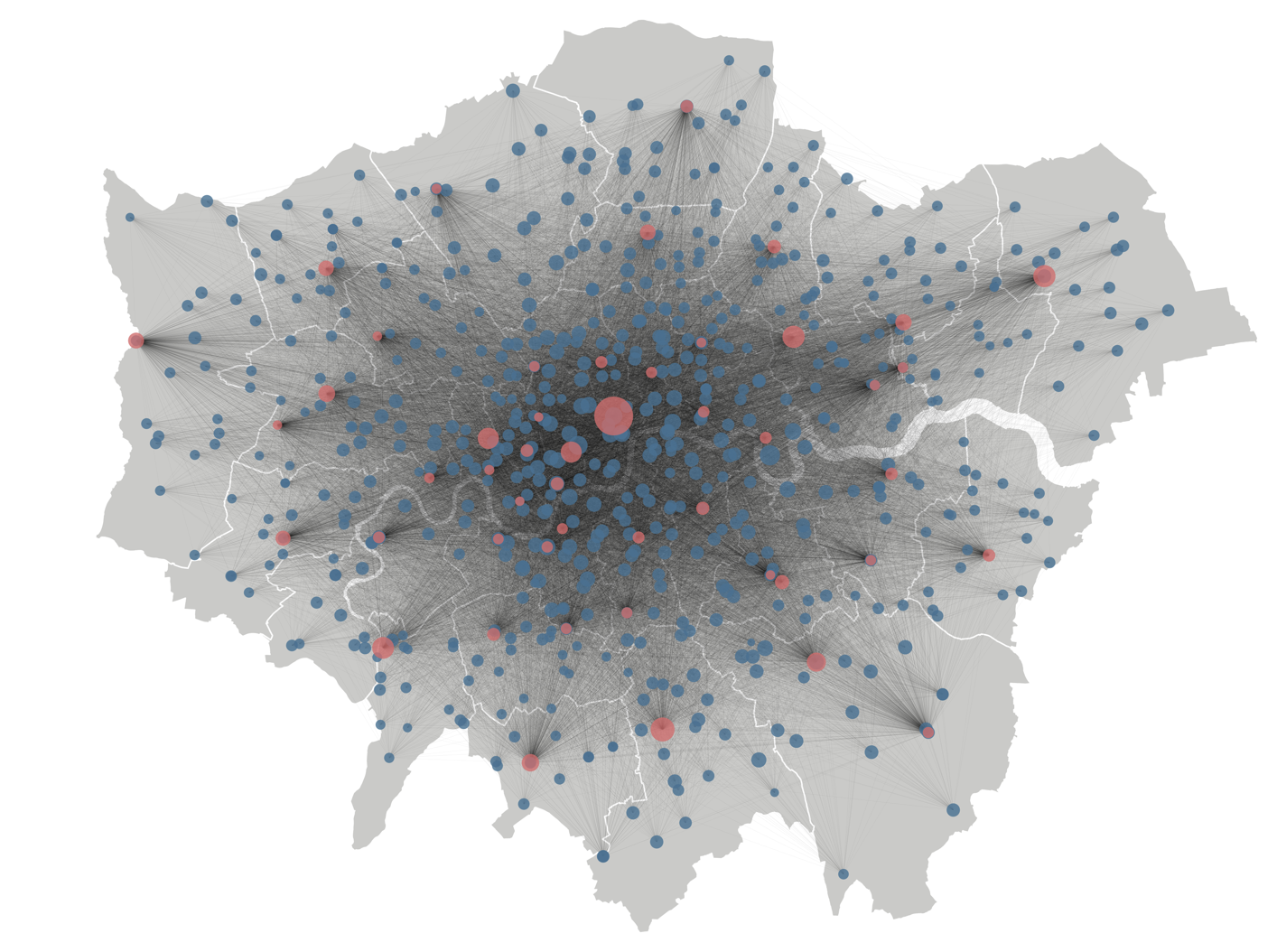}
	\end{minipage}
	\hspace{.0266\textwidth}
	\begin{minipage}[!t]{.333\textwidth}
		\includegraphics[width=\textwidth]{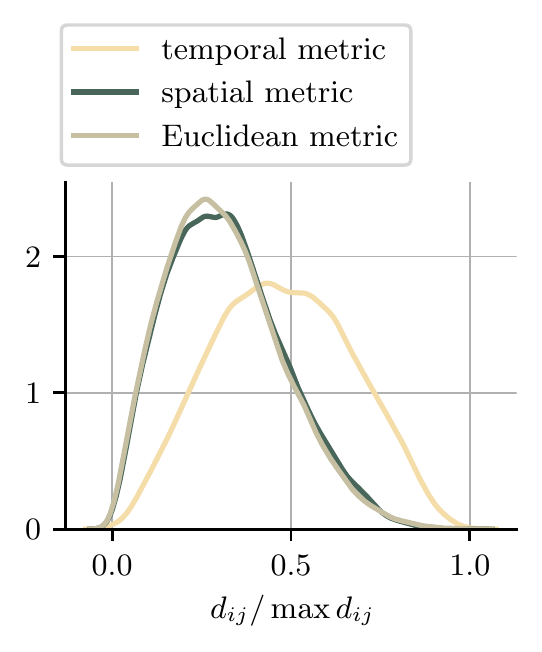}
	\end{minipage}
	\caption{Left: Visualisation of the London dataset, as described in the text. Blue dots are the origin zone sizes for all $N=625$ wards, red dots are retail floor space sizes for $M=49$ retail centres. The network edge weights are given by eq. [\ref{eq:convenience_factor}], using the minimum transport time between locations for public transport or driving. Background map of London boroughs: \cite{GLA_GIS}. Right: the densities of the distance distributions for three different metrics: temporal, spatial, and Euclidean.}
	\label{fig:London_dataset}
\end{figure*}
We apply our methodology to a dataset of economic activity in Greater London and compare the results to a study of the same dataset, performed using an MCMC approach \cite{Ellam_2018}. The data consists of origin zone sizes $\mathbf{O}$, destination zone sizes $\mathbf{W}$, and two different convenience matrices $\mathbf{C}$ for the region for the period around 2015/2016.

\paragraph*{Origin zones:} The origin zone sizes are given by the total spending budget of each of the $N=625$ electoral wards in Greater London as classified by the Greater London Authority (GLA) \cite{GLA_ward_data}. Ward-level income is calculated as the product of the number of households times the average household income; however, households do not spend their entire budget on retail and service goods: according to the Office for National Statistics \cite{ONS_household_budget}, between 2015--2017 London households on average spent between 21\%--30\% of their total budget on comparison, service, and convenience goods\footnote{The 30\% figure comprises food and non-alcoholic drinks, alcoholic drinks and tobacco, clothing and footwear, household goods and services, and miscellaneous goods and services; this figure drops to 21\% when excluding food and non-alcoholic drinks. See \cite{ONS_household_budget}, tab. A33.}, hence we multiply the income figures by 0.21 to obtain the origin size (we are excluding smaller retail zones from the dataset, see below, and hence choose to exclude food costs from the budget figure, see footnote). For example, for the City of London, 6680 households and an average household income of \textsterling 63.620/a result in a total ward-level income of about \textsterling425~million~per~annum, and hence an origin size of \textsterling89 million/annum. In order to prevent numerical overflow, we use units of \textsterling $10^8$/a for the origin zones. Recall that scaling the origin zone values only affects the resulting prediction for $\kappa$, as $\alpha$ and $\beta$ are unitless. It should be noted that while the population figures are for 2015, the income figures are only for 2012/2013. 

\paragraph*{Destination zones:} The destination zone sizes are the total occupied retail floor space sizes in $\mathrm{m}^2$ for all $M=49$ town centres classified as either `international', `metropolitan', or `major' by the Greater London Authority (GLA) \cite{ GLA_health_check_report}. In the GLA report \cite{GLA_health_check_report} this retail floor space comprises comparison, convenience, and service retail (as given by \cite{GLA_health_check_report}, tab. 1.1.1a--1.1.1c); for example, for the West End we obtain a total occupied retail floor space in 2016 of 474.456 $\mathrm{m}^2$. In order to prevent numerical overflow (resulting from large values of $W_j^\alpha$ potentially occurring during the training), we use units of $10^5 \mathrm{m}^2$ for the destination zone sizes. With this choice of units for the origin and destination zone sizes, $\kappa$ will be given in units of \textsterling1000/a. The true value resulting from this data is $\kappa = \sum_i O_i / \sum_k W_k = $\textsterling8301/a.
\begin{figure*}[!ht]
\begin{minipage}{0.3\textwidth}
    {\textbf{(a)} Low training noise  ($\sigma = 0.014$)}
\end{minipage}
\begin{minipage}{0.7\textwidth}
    \includegraphics[width=\textwidth]{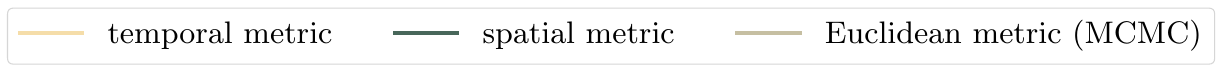}
\end{minipage}
\begin{minipage}{\textwidth}
    \includegraphics[width=0.332\textwidth]{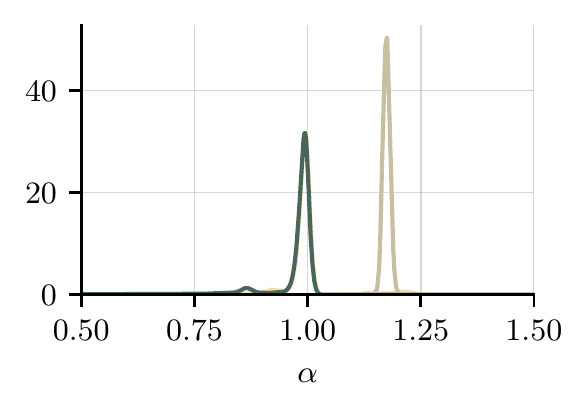}
    \includegraphics[width=0.332\textwidth]{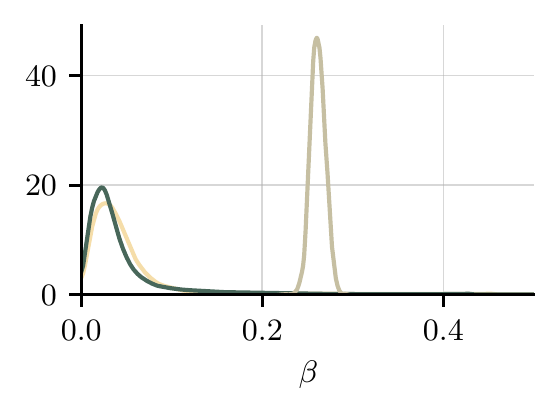}
    \includegraphics[width=0.332\textwidth]{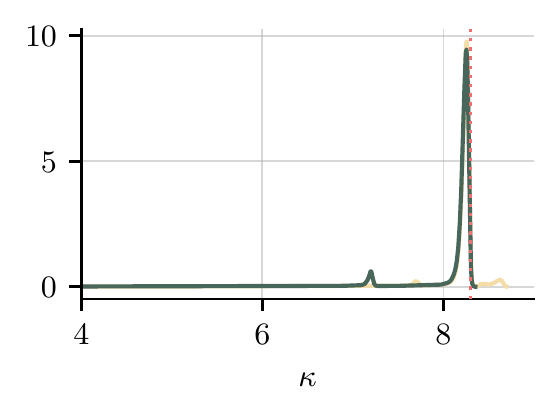}
\end{minipage}

\begin{minipage}{\textwidth}
    {\textbf{(b)} High training noise ($\sigma = 0.14$)}
\end{minipage}
\begin{minipage}{\textwidth}
    \includegraphics[width=0.332\textwidth]{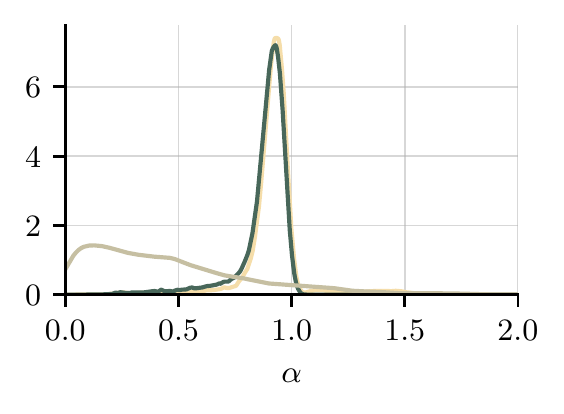}
    \includegraphics[width=0.332\textwidth]{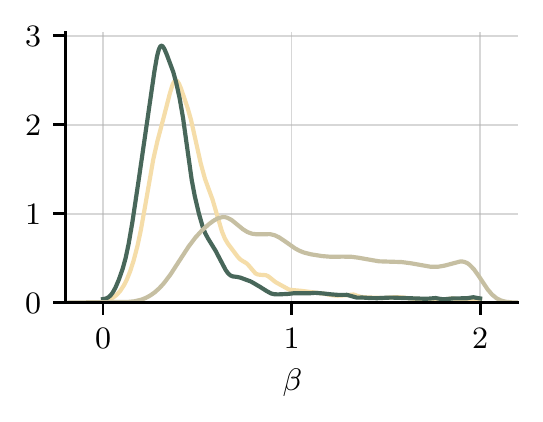}
    \includegraphics[width=0.332\textwidth]{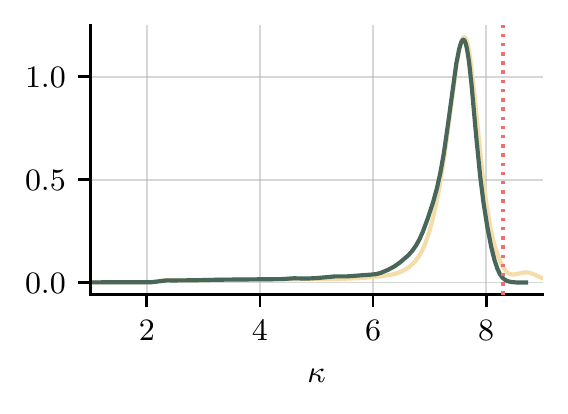}
\end{minipage}
   
\caption{Marginal densities for $\boldsymbol{\lambda} = (\alpha, \beta, \kappa)$ for different network metrics and two noise regimes. Densities are obtained by running the model from 20 different initial seeds for 10000 epochs each, and smoothed with a Gaussian kernel. The true value of $\kappa$ is given by the red dotted line. Runtime for each metric and noise regime: 1 min 8 secs. Also shown are the densities for $\alpha$ and $\beta$ obtained by running the MCMC analysis presented in \cite{Ellam_2018} (light brown). Runtime: 3 hrs 40 minutes (low noise), 7 hrs 20 minutes (high noise).}
\label{fig:London_parameter_estimates}
\end{figure*}
\begin{table*}[t]
\vspace{3mm}
\centering
\begin{tabular}{ll|ccc}
Noise regime &   & Temporal Metric & Spatial Metric & Euclidean metric (MCMC) \\
\midrule
\multirow{3}{*}{$\sigma=0.014$ (low noise)} & predictions $(\hat{\alpha}, \hat{\beta}, \hat{\kappa})$ & $(0.99, 0.02, 8.26)$ & $(0.99, 0.02, 8.25)$ & $(1.18, 0.28, (8.3))$ \\
& expected MSPE & $(1.7 \pm 1.2) \times 10^{-8}$ & $(1.6 \pm 1.1) \times 10^{-8}$ & $(4.5 \pm 0.1) \times 10^{-5}$ \\ 
 & compute time & 1 min 8 secs & 1 min 8 secs & 3 hrs 40 mins \\
\midrule
\multirow{3}{*}{$\sigma=0.14$ (high noise)} & predictions $(\hat{\alpha}, \hat{\beta}, \hat{\kappa})$ & $(0.92, 0.54, 7.45)$ & $(0.92, 0.39, 7.39)$ & $(0.12, 0.75, (8.3))$ \\
& expected MSPE & $(1.8 \pm 1.3) \times 10^{-6}$ & $(1.7 \pm 1.1) \times 10^{-6}$ & $(2.3 \pm 0.2)\times 10^{-4}$ \\ 
 & compute time & 1 min 8 secs & 1 min 8 secs & 7 hrs 20 mins \\
\bottomrule
\end{tabular}
\caption{Comparison of predicted values, expected calibration MSPE, and compute times}
\label{tab:HW_results}
\end{table*}
\paragraph*{Cost network:} For the cost network $\mathbf{C}$ we use the Google Distance Matrix API\footnote{\href{https://developers.google.com/maps/documentation/distance-matrix}{developers.google.com/maps/documentation/distance-matrix}} to extract travel data from Google Maps. The API can be used to extract both travel distances and travel times for large batches of search queries for all the transport modes available on Google Maps (driving, public transport, walking, cycling). We only consider two modes of transport: driving and public transport. The API does not provide past data, so travel time and distance data reflect the state of the network in June 2022. In order to restore some level of comparability, we set the travel date to a Sunday, in order to blend out the newly opened Elizabeth line (which in June 2022 did not yet run on Sundays). We define the distance $d_{ij}$ between two nodes as the \emph{shorter} of the two transport modes considered: for example, travelling from Kentish Town to the West End takes about 19.7 minutes by public transport, and 22.35 minutes by car, so we set the distance to be 19.7 minutes (of course, this assumes that everyone can choose between both modes of transport, and does not factor in the added cost of driving in London.) Note that the Google Maps public transport mode defaults to walking for very short distances. We consider two different metrics: a temporal metric and a spatial metric. The convenience factors are derived in \cite{Wilson_1967} as 
\begin{equation}
	c_{ij} = e^{- d_{ij}/\tau}
	\label{eq:convenience_factor}
\end{equation}
where $d_{ij}$ is the distance in the metric in question and $\tau = \sup_{i,j} d_{ij}$ the time/length scale, ensuring a unitless exponent. The network is then given by $\mathbf{C} = (c_{ij})$. Figure \ref{fig:London_dataset} visualises the dataset: the red nodes are the destination zones $W_j$, the blue dots are the origin zones $O_i$, and the network edge widths represent the cost matrix, using the temporal metric. The large, central destination zone is the London West End, by far the largest retail zone in London. On the right, the distance distributions for the two metrics are shown: it is interesting to observe that the temporal metric has a distribution more or less symmetrically centred around 0.5, whereas the distance metric is much more heavily skewed to the left. In practice this means that zones are statistically closer spatially than temporally. This may be an artefact of our choice of transportation mode: since users must use either public transport or drive, zones within realistic walking or cycling distance appear further away than they actually are. This is substantiated by the fact that more zones lie within zero spatial distance of each other than within zero temporal distance. We expect the choice of metric to only affect the predictions for $\beta$, the convenience parameter coupling the dynamics to the network.
\paragraph*{Results} We estimate the marginal densities for the parameters for both metrics (temporal and spatial), and compare the results to those obtained in \cite{Ellam_2018}. In that study, the authors considered two underlying noise levels in the data, $\sigma=\sqrt{2} \times 10^{-2}$ and $\sigma=\sqrt{2} \times 10^{-1}$, and estimated densities for $\alpha$ and $\beta$. We produce estimates for $\alpha$, $\beta$, and $\kappa$, using the same values for the noise in the forward pass of the scheme. The resulting marginals are shown in fig. \ref{fig:London_parameter_estimates}, along with the densities obtained from running the MCMC scheme from \cite{Ellam_2018}. We train our model for 10.000 epochs over 20 different initialisations to obtain the same number of samples. It should be noted that the authors in \cite{Ellam_2018} used Euclidean distances for their study, which however in fig. \ref{fig:London_parameter_estimates} can be seen to have much the same distribution as our spatial metric. 

Table \ref{tab:HW_results} summarises the results. While the MCMC scheme produced very different estimates depending on the choice of $\sigma$, our method shows good consistency across noise regimes. As is to be expected, the predictions are independent of the choice of metric for $\alpha$ and $\kappa$. We obtain good to fair estimates for $\kappa$, as well as a consistent maximum likelihood estimate for $\alpha$. We assess the quality of the predictions by using the most likely values to calibrate the numerical solver eq. [\ref{eq:HW_model}], and computing the expected mean squared prediction error over 100 runs (thereby accounting for the random noise involved). Since \cite{Ellam_2018} did not estimate $\kappa$, we use its true value for the calibration run with the MCMC values. As can be seen, our method improves on the MCMC prediction by three orders of magnitude in the low noise, and two orders of magnitude in the high noise regime. Note also that our method produces density estimates in just over one minute, while the MCMC scheme took between 3.7 hours and 7.3 hours to run; this represents a performance improvement of two orders of magnitude. The speed-up is not attributable to parallelisation: in CPU time, our method took 9 minutes to run, still representing two orders of magnitude performance increase.

\section{Conclusions}
In this article, we outlined an approach to estimating probability densities for parameters of differential equations using neural nets. We considered a broad spectrum of datasets, including time series data of different lengths, steady state data with only a single time frame, noiseless and noisy data, synthetically generated, and real data; the method was applied to a non-convex problem with two global minima of the loss function (Harris-Wilson model), as well as to a situation with data-model mismatch (SIR model). In all cases, the neural net quickly and reliably found densities for the relevant parameters. We assessed the quality of the predictions by comparing the parameter predictions to their true values, in the case of the Harris-Wilson model, comparing our model's prediction errors to that of a study using classical MCMC techniques. Our model significantly outperformed that study in terms of the accuracy of the predicted data, while being computationally much faster to run: in both regards, improvements covered two to three orders of magnitude. The use of numerical solvers for the forward pass allows giving estimates for different levels of assumed noise in the data, which can help to obtain realistic confidence bounds on predictions. However, in the SIR study we also showed the ability of our method to learn the noise level itself. We demonstrated the importance of (1) analysing a model's dynamical range before training, (2) ensuring numerical stability throughout the training process, and (3) adapting the neural net architecture to the problem. For the Harris-Wilson model, this was achieved by (1) only training on data from the dynamical range $0 < \nu < 1$, (2) scaling the origin and destination zone sizes appropriately to prevent numerical overflow, and (3) ensuring positive predictions by using the modulus as the activation function.

Due to the method's performance and relative mathematical simplicity, we hope that it will prove a powerful tool across the quantitative sciences, such as in quantitative sociology or computational epidemiology. Our work opens up fruitful lines for further research: in particular, we have not yet considered how the model performance scales with the size of $\boldsymbol{\lambda}$. How, for instance, will the model fare when predicting very large parameter sets, such as graph adjacency matrices? This will be the subject of further investigation by the authors.
\vspace{0.5cm}
\hrule 
\vspace{0.5cm}
\small{\noindent {\bf Acknowledgements} TG was funded by the University of Cambridge School of Physical Sciences VC Award via DAMTP and the Department of Engineering, and supported by EPSRC grant EP/P020720/2. The work of GP was partially funded by EPSRC grant EP/P031587/1, and by J.P. Morgan Chase \& Co through a Faculty Research Award 2019 and 2021. MG was supported by EPSRC grants EP/T000414/1, EP/R018413/2, EP/P020720/2, EP/R034710/1, \hspace{\fill}\linebreak EP/R004889/1, and a Royal Academy of Engineering Research Chair. The authors would also like to thank Dr. Claudia Totzeck (University of Wuppertal) for useful discussions, whose paper \cite{Goettlich_2021} was a starting point for this work, as well as Louis Ellam for his assistance in procuring the GLA dataset.}

\vspace{0.5cm}
\hrule

% % Bibliography
\bibliographystyle{bibstyle}
\bibliography{library.bib}

\appendix
\onecolumn
\cleardoublepage\phantomsection\addcontentsline{toc}{section}{Supporting Information}

\renewcommand\thesection{S\arabic{section}}
\setcounter{section}{0}
{\sffamily \noindent \huge \bfseries Supporting Information}
\vspace{2cm}
\renewcommand\thefigure{S\arabic{figure}}
\setcounter{figure}{0}
\section*{Methodology}

\lstset{frame=l,
  backgroundcolor=\color{gray},
  aboveskip=3mm,
  belowskip=3mm,
  showstringspaces=false,
  columns=flexible,
  basicstyle={\small\ttfamily},
  numbers=none,
  numberstyle=\tiny\color{gray},
  keywordstyle=\color{blue},
  commentstyle=\color{dkgreen},
  stringstyle=\color{mauve},
  breaklines=true,
  breakatwhitespace=true,
  tabsize=3
}
\lstdefinestyle{yaml}{
     basicstyle=\color{blue}\footnotesize,
     rulecolor=\color{black},
     string=[s]{'}{'},
     stringstyle=\color{blue},
     comment=[l]{:},
     commentstyle=\color{black},
     morecomment=[l]{-}
 }

\subsection*{Neural Networks: Notation and Terminology}
A \emph{neural network} is a sequence of length $L \geq 1$ of concatenated transformations. Each \emph{layer} of the net consists of $L_i$ \emph{neurons}, connected through a sequence of \emph{weight matrices} $\mathbf{W}_i~\in~\mathbb{R}^{L_{i+1} \times L_i}$. Each layer applies the transformation
\begin{equation*}
	\sigma_i(\mathbf{W}_i \mathbf{x} + \mathbf{b}_i)
\end{equation*}
to the input $\mathbf{x}$ from the previous layer, where $\mathbf{b}_i \in \mathbb{R}^{L_{i+1}}$ is the \emph{bias} of the $i$-th layer. The function $\sigma_i: \mathbb{R}^{L_{i+1}} \rightarrow \mathbb{R}^{L_{i+1}}$ is the \emph{activation function}; popular choices include the \emph{rectified linear unit} (ReLU) activation function $\sigma(x) = \max(x, 0)$, and the \emph{sigmoid activation function} $\sigma(x) = (1 + e^{-x})^{-1}$. A neural net has an \emph{input layer}, an \emph{output layer}, and \emph{hidden layers}, which are the layers in between the in- and output layers. If a network only has one hidden layer, we call it \emph{shallow}, else we call the neural net \emph{deep}.

\begin{figure}[h]
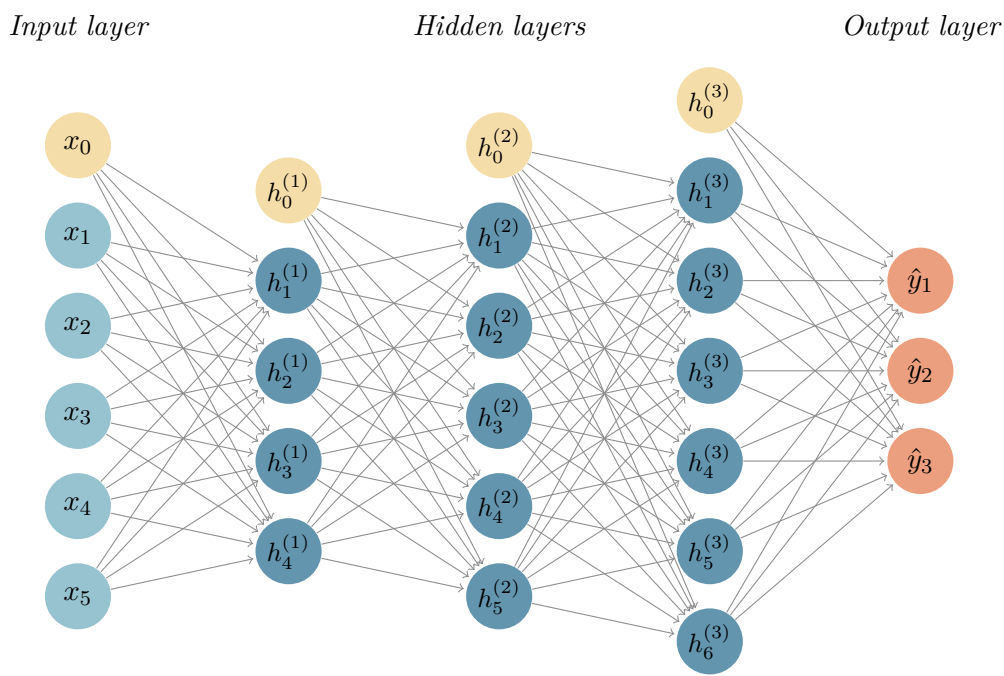

\centering
\definecolor{yellow}{HTML}{F5DDA9}
\definecolor{green}{HTML}{48675A}
\definecolor{red}{HTML}{ec7070}
\definecolor{blue}{HTML}{2F7194}
\definecolor{lightblue}{HTML}{97c3d0}
\definecolor{orange}{HTML}{EC9F7E}
\definecolor{brown}{HTML}{C6BFA2}
 \begin{neuralnetwork}[height=7, nodesize=25, 
 		nodespacing=12mm, layerspacing=28mm]
  		
  		\tikzstyle{bias neuron}=[neuron, fill=yellow!100];
  		\tikzstyle{hidden neuron}=[neuron, fill=blue!75];
  		\tikzstyle{output neuron}=[neuron, fill=orange!100];
  		\tikzstyle{input neuron}=[neuron, fill=lightblue!100];
        
        \newcommand{\x}[2]{$x_#2$}
        \newcommand{\y}[2]{$\hat{y}_#2$}
        \newcommand{\hfirst}[2]{\small $h^{(1)}_#2$}
        \newcommand{\hsecond}[2]{\small $h^{(2)}_#2$}
        \newcommand{\hthird}[2]{\small $h^{(3)}_#2$}
        
        \inputlayer[count=5, bias=true, title={\it Input layer}, text=\x]
        \hiddenlayer[count=4, bias=true, text=\hfirst] \linklayers
        \hiddenlayer[count=5, bias=true, title={\it Hidden layers}, text=\hsecond] \linklayers
        \hiddenlayer[count=6, bias=true, text=\hthird] \linklayers
        \outputlayer[count=3, title={\it Output layer}, text=\y] \linklayers
\end{neuralnetwork}
\vspace{5mm}
\caption{Example of a deep neural network with 3 hidden layers. The inputs (light blue nodes) are passed through the layers, with links between layers representing the weight matrices $\mathbf{W}$. Each layer also applies a \emph{bias} (yellow nodes), with the network finally producing an output (orange).}
\label{diag:neural_net}
\end{figure}

\subsection*{Details on the code}
The code is uploaded to the \href{https://github.com/ThGaskin/NeuralABM}{Github repository} as given in the main text.
\subsubsection*{Installation}
Detailed installation instructions are given in the repository. First, clone the repository, install the \href{https://utopya.readthedocs.io/en/latest/}{\texttt{utopya}} package and all the required additional components into a virtual environment, for example via \texttt{PyPi}. In particular, install \href{https://pytorch.org}{\texttt{pytorch}}. Enter the virtual environment. Then, from within the project folder, register the project:
\begin{lstlisting}
utopya projects register .
\end{lstlisting}
You should get a positive response from the utopya CLI and your project should appear in the project list when calling:
\begin{lstlisting}
utopya projects ls
\end{lstlisting}
Note that any changes to the project info file need to be communicated to utopya by calling the registration command anew. You will then have to additionally pass the \texttt{--exists-action overwrite} flag, because a project of that name already exists. See 
\begin{lstlisting}
utopya projects register --help
\end{lstlisting}
for more information.
Finally, register the SIR model via
\begin{lstlisting}
utopya models register from-manifest models/SIR/SIR_info.yml
\end{lstlisting}
(and the Harris-Wilson model accordingly).
\subsubsection*{Running the code}
To run the code, execute the following command:
\begin{lstlisting}
utopya run <model_name>
\end{lstlisting}
By default, this runs the model with the settings in the \texttt{<model\_name>\_cfg.yml} file. All data and the plots are written to an output directory, typically located in \texttt{\textasciitilde/utopya\_output}. To run the model with different settings, create a \texttt{run\_cfg.yml} file and pass it to the model like this:
\begin{lstlisting}
utopya run <model_name> path/to/run_cfg.yml
\end{lstlisting}
This is recommended rather than changing the default settings, because the defaults are parameters that are known to work and you may wish to fall back on in the future. 

Plots are generated using the plots specified in the \texttt{<model\_name>\_plots.yml} file. These too can be updated by creating a custom plot configuration, and running the model like this:
\begin{lstlisting}
utopya run <model_name> path/to/run_cfg.yml --plots-cfg path/to/plot_cfg.yml
\end{lstlisting}
See the \href{https://docs.utopia-project.org/html/getting_started/tutorial.html}{Utopia tutorial} for more detailed instructions.

All the images in this article can be generated using so-called \emph{configuration sets}, which are complete bundles of both run configurations and evaluation configurations. For example, to generate the four frames of the SIR model, you can call
\begin{lstlisting}
utopya run SIR --cfg-set ABM_data
\end{lstlisting}
This will run and evaluate the \texttt{SIR} model with all the settings from the \texttt{SIR/cfgs/ABM\_data/run.yml} and \texttt{eval.yml} configurations.
\subsubsection*{Parameter sweeps}
Parameter sweeps are automatically parallelised by \texttt{utopya}, meaning simulation runs are always run concurrently whenever possible. The data is automatically stored and loaded into a data tree. To run a sweep, simply add a \texttt{!sweep} tag to the parameters you wish to sweep over, and specify the values, along with a default value to be used if no sweep is performed:
\begin{lstlisting}
param: !sweep 
  default: 0
  values: [0, 1, 2, 3]
\end{lstlisting}
Then in the run configuration, add the following entry:
\begin{lstlisting}
[bgcolor=gray]{yaml}
perform_sweep: true
\end{lstlisting}
Alternatively, call the model with the flag \texttt{-{}-run-mode sweep}. The configuration sets used in this work automatically run sweeps whenever needed, so no adjustment is needed to recreate the plots used in this work.
\subsubsection*{Initialising the neural net}
The neural net is controlled from the \texttt{NeuralNet} entry of the configuration:
\begin{lstlisting}
NeuralNet:
  num_layers: 4
  nodes_per_layer: 
    default: 20
    layer_specific:
      1: 10
      2: 15
  biases:                   # optional; if this entry is omitted no biases are used
    default: ~              # default is None (indicted by a tilde in YAML)
    layer_specific:
      0: default            # use pytorch default (xavier uniform)
      -1: [-1, 1]           # uniform initialisation on a custom interval
  activation_funcs:         
    default: sigmoid   
    layer_specific:         # optional
      1: tanh
      3:
        name: HardTanh     # you can also pass a function that takes additional args and/or kwargs
        args:
          - -2  # min_value
          - +2  # max_value
  learning_rate: 0.001      # optional; default is 0.001
  optimizer: SGD            # optional; default is Adam           
\end{lstlisting}
\texttt{num\_layers} specifies the depth of the net; \texttt{nodes\_per\_layer} controls the architecture: provide a \texttt{default} size, and optionally any deviations from the default under \texttt{layer\_specific}. The keys of the \texttt{layer\_specific} entry should be indices of the layer in question. The optional \texttt{biases} entry determines whether or not biases are to be included in the architecture, and if so how to initialise them. A default and layer-specific values can again be passed. Setting an entry to \texttt{default} initialises the values using the pytorch default initialiser, a Xavier uniform initialisation. Passing a custom interval instead initialises the biases uniformly at random on that interval, and passing a tilde \texttt{\textasciitilde} (\texttt{None} in YAML) turns the bias off. \texttt{activation\_funcs} is a dictionary specifying the activation functions on each layer, following the same logic as above. Any \href{https://pytorch.org/docs/stable/nn.html}{pytorch activation function} is permissible. If a function requires additional arguments, these can be passed as in the example above. Lastly, the \texttt{optimizer} keyword takes any argument \href{https://pytorch.org/docs/stable/optim.html#base-class}{allowed in pytorch}. The default optimizer is the Adam optimizer \cite{Kingma_2014} with a learning rate of $0.001$.

The neural net can be initialised from different initial values in the parameter space by changing the \emph{random seed} in the configuration:
\begin{lstlisting}
seed: 42
\end{lstlisting}
Sweeping over different initialisations is achieved by sweeping over the seed, as described in the previous section.
\subsubsection*{Training the neural net}
The \texttt{Training} entry of the configuration controls the training process:
\begin{lstlisting}
Training:
  batch_size: 2  
  to_learn: [param1, param2]
  true_parameters:
    param3: 0.4
  loss_function:
    name: MSEloss
    # can pass additional args and kwargs here ...
\end{lstlisting}   
You must specify which parameters to learn, and pass the true values for the others if you are not learning all parameters. Under the \texttt{loss\_function} key you can specify the loss function to use, and pass any arguments or keyword arguments it may require using an \texttt{args} or \texttt{kwargs} key. You can use any available \href{https://pytorch.org/docs/stable/nn.html#loss-functions}{pytorch loss function}. 

\section*{Application to time series data: diffusive SIR model
of epidemics}
\subsection*{Neural Network Architecture}
For this section, we choose a shallow neural net with 20 neurons in the hidden layer and the absolute value function as the activation function to ascertain parameter predictions are positive. We initialise the biases uniformly at random with values in the unit interval $[0, 1]$ to further push the neural network towards positive predictions. We use the Adam optimizer with a learning rate of $0.001$. The corresponding configuration looks like this:

\begin{lstlisting}
SIR:
  NeuralNet:
    num_layers: 1
    nodes_per_layer: 
      default: 20
    biases: 
      default: [0, 1]  # Initialise biases on [0, 1]
    activation_funcs:
      default: linear                 
      layer_specific:
        -1: abs # Modulus activation func
    learning_rate: 0.002      
    optimizer: Adam           
  Training:
    batch_size: 90            # the time series has length 100
    to_learn: [p_infect, t_infectious, sigma]
    loss_function:
      name: MSELoss
\end{lstlisting}
 
\subsection*{Training}
We train the model using the following operation: let $\boldsymbol{\varphi}(t) = (S(t), I(t), R(t))$ be the current state of the model; then in each iteration, we do
\begin{equation}
    \boldsymbol{\varphi}(t+1) = \mathrm{ReLU}(\boldsymbol{\varphi}(t) + \mathrm{d}\boldsymbol{\varphi}(t)),
\end{equation}
where $\mathrm{d}\boldsymbol{\varphi}(t)$ is given by 
\begin{equation}
    \mathrm{d}\boldsymbol{\varphi}(t) =
    \begin{pmatrix}
    -\hat{\beta} & -\hat{w} \\ 
    \hat{\beta} & \hat{w} -  f(\gamma\hat{\tau}, t) \\
    0 & f(\gamma\hat{\tau}, t)
    \end{pmatrix}
    \begin{pmatrix}
    S(t)I(t) \\ I(t)
    \end{pmatrix}
\end{equation}
The ReLU function ensures that densities do not become negative.
The estimated value $\hat{w}$ is given by $\hat{\sigma} X$, where $\hat{\sigma}$ is the neural net prediction for the noise, and $X \sim N(0, 0.1)$ a normally distributed random variable with variance $0.1$. The function $f$ is given by
\begin{equation}
    f(s, t) = \dfrac{1}{s} \mathrm{sigmoid}(k t/s)
\end{equation}
with $k = 1000$. It approximates a step function, ensuring that recovery only begins after a certain period. $\gamma$ is a scaLling factor, designed to ensure the three estimated parameters are roughly of the same order of magnitude. This increases the speed of the neural net's convergence to a loss function minimum, as the parameters are closer together, but is not required to achieve reliable results. When writing the data, $\hat{\tau}$ is scaled back to its original dimension. We choose $\gamma=10$. 
\subsection*{Running the code}
The training data for the ABM is provided in the \texttt{data/SIR} folder, and is provided in \texttt{hdf5} format. To train the neural net, simply run the following command (make sure you are in the project folder, otherwise change the path to the dataset to an absolute path):
\begin{lstlisting}
utopya run SIR --cfg-set Predictions
\end{lstlisting}
This will load the ABM data and run the model 20 times from different initialisations, finally producing the plots from the manuscript. You can change the number of sweeps, and the initial seeds, by changing the \texttt{seed} entry in the \texttt{run.yml} configuration. By default, it looks like this:
\begin{lstlisting}
seed: !sweep
  default: 0
  range: [20]
\end{lstlisting}

\section*{Application to a non-convex problem: the Harris-Wilson model of economic activity}

\subsection*{Neural Network Architecture}
The neural network architecture is the same as for the SIR section, using the same optimizer with the same learning rate of $0.002$. We initialise the biases of the neural net in the interval $[0, 4]$:
\begin{lstlisting}
HarrisWilson:
  NeuralNet:
    num_layers: 1
    nodes_per_layer: 
      default: 20
    activation_funcs:
      default: linear
      layer_specific:
        -1: abs
    biases:
      default: [0, 4]
    learning_rate: 0.002
    optimizer: Adam
  Training:
    to_learn: [alpha, beta, kappa]
    batch_size: 1
    true_parameters:
      sigma: 0
    loss_function:
      name: MSELoss
\end{lstlisting}
In the noiseless case we set the training noise to $0$, as we are not learning the noise; for the noisy runs, we also learn the noise.
\subsection*{Training}
We use the following matrix form of the Harris-Wilson equations to train the neural net. Let $\mathbf{D} \in \mathbb{R}^M$, $\mathbf{O} \in \mathbb{R}^N$, $\mathbf{W} \in \mathbb{R}^M$ be the demand vector, origin zone size vector, and destination zone size vector respectively. The the dynamics are given by
\begin{equation}
	\mathbf{D} = \mathbf{W}^\alpha \odot \left[\left(\mathbf{C}^\beta\right)^T (\mathbf{O} \odot \mathbf{Z})  \right] \in \mathbb{R}^M,
\end{equation}
with $\odot$ indicating the Hadamard product, and elementwise exponentiation. $\mathbf{Z} \in \mathbb{R}^N$ is the vector of normalisation constants
\begin{equation}
	\mathbf{Z}^{-1} =  \mathbf{C}^\beta \mathbf{W}^\alpha.
\end{equation} 
The dynamics then read
\begin{equation}
	\dot{\mathbf{W}} = \mathbf{W} \odot \epsilon(\mathbf{D} - \kappa \mathbf{W})
	\label{eq:matrix_formulation}
\end{equation}
with given initial conditions $\mathbf{W}(t=0) = \mathbf{W}_0$.
This formulation is more conducive to machine learning purposes, since it contains easily differentiable matrix operations and does not use for-loop iteration.
\subsection*{Performance analysis}
\begin{figure}[t]
	\includegraphics[width=\textwidth]{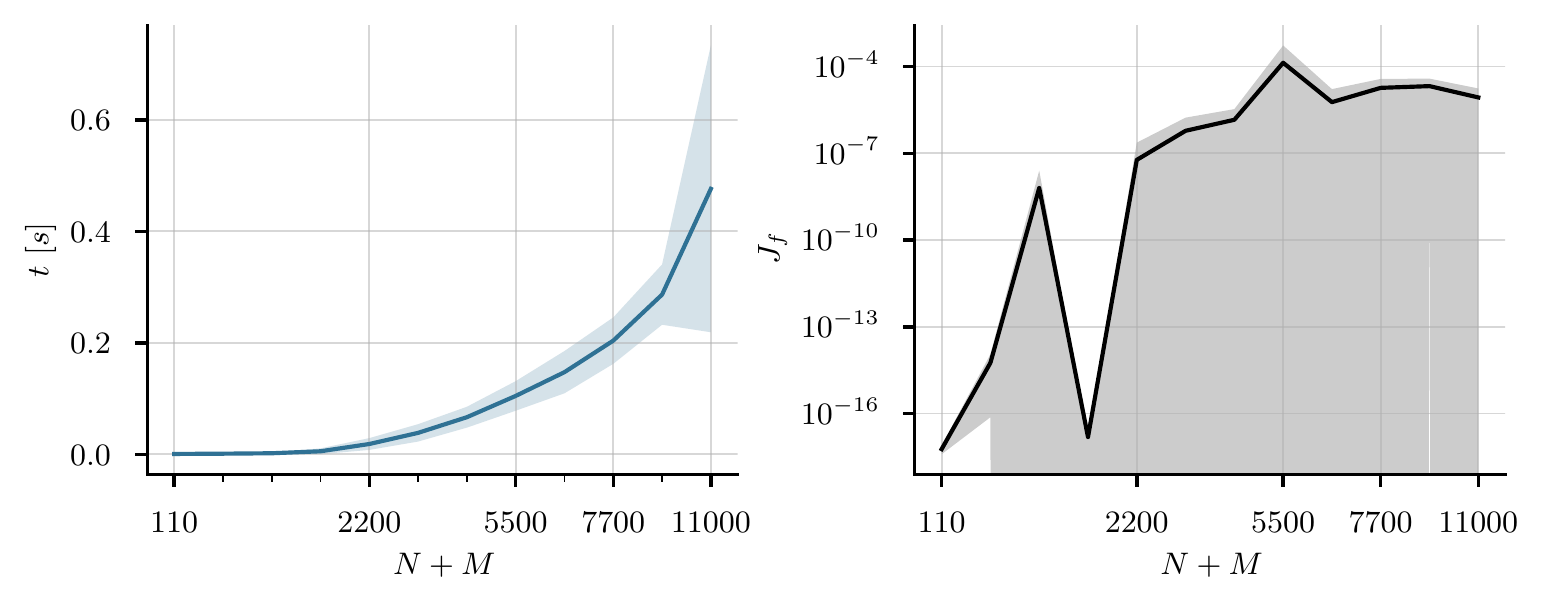}
	\caption{Performance analysis of the model as a function of the network size $N+M$. Left: time to complete a single epoch ($L=1$, $B=1$). Right: training loss after 6000 epochs.}
	\label{fig:HW_performance_analysis}
\end{figure}

Figure \ref{fig:HW_performance_analysis} shows a performance analysis of the model as a function of the network size $N+M$. On the left, the time a neural net with 20 neurons takes to run a single epoch ($L=1$, $B=1$). Each data point is an average over $10$ different initialisations, each with 6000 epochs. On the right, the loss after for each size after 6000 epochs is shown, averaged over 10 different initialisations. The shaded area represents the standard deviation. \textcolor{blue}{Training was performed on the CPU of a standard laptop.}
\subsection*{Marginal densities}
\begin{figure*}[h!]
    \begin{minipage}{\textwidth}
    \centering
    \includegraphics[width=0.4\textwidth]{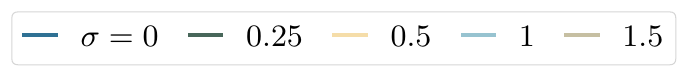}
    \end{minipage}
    \includegraphics[width=0.32\textwidth]{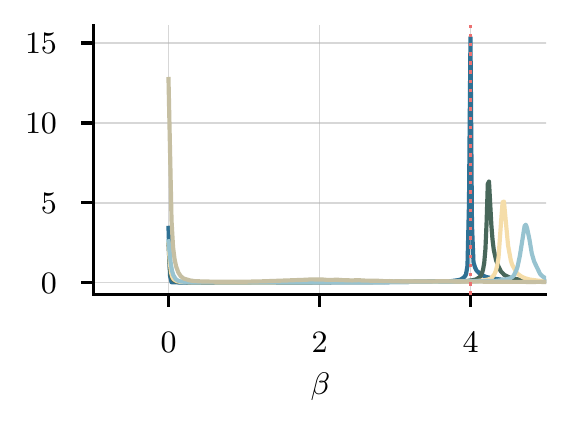}
     \includegraphics[width=0.32\textwidth]{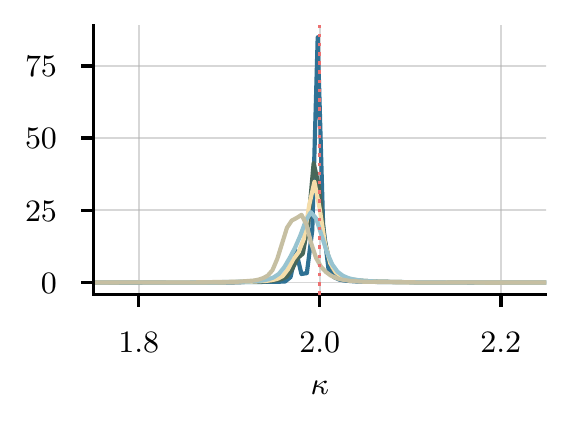}
    \includegraphics[width=0.32\textwidth]{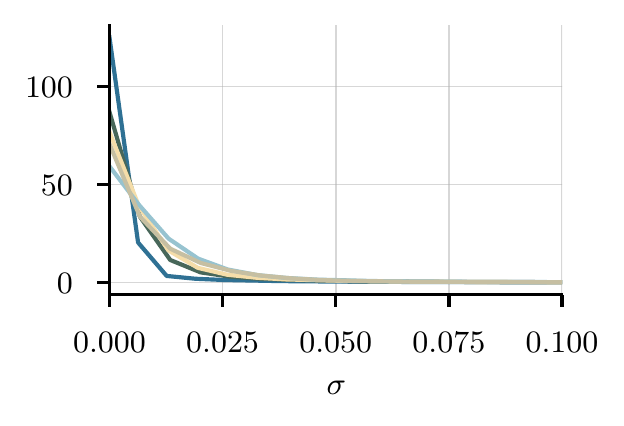}
	\caption{Marginal densities for $\beta$, $\kappa$, and $\sigma$ for different levels of noise in the training data, smoothed with a Gaussian kernel. All configurations as in figure \ref{fig:HW_synthetic_marginals}. Red dotted lines: true values.}
	\label{fig:HW_additional_marginals}
\end{figure*}

Figure \ref{fig:HW_additional_marginals} shows the corresponding marginals for $\beta$, $\kappa$, and the noise level $\sigma$ for figure 7 in the main manuscript.
\subsection*{Marginal densities for $\beta$, $\kappa$, $\sigma$}
\subsection*{Datasets}
The \texttt{data/HarrisWilson} folder contains several datasets, both real and synthetic, that can be used to train the neural net and learn 
parameters for the Harris-Wilson equations. Simply set the \texttt{load\_from\_dir} key in your \texttt{run.yml} file to point
to the folder containing the data: the data will automatically be loaded and 
the model trained on that data. 

\subsubsection*{Synthetic data}
The \texttt{synthetic\_data} folder contains synthetically generated networks, origin sizes, and destination sizes, both with and without noise. The name of the 
folder indicates the network size, i.e. \texttt{N\_100\_M\_10} means $N=100$ and $M=10$. However, all folders also include a 
\texttt{config.yml} file detailing the specific configurations for each dataset.

\subsubsection*{London dataset}
The \texttt{London\_data} folder contains datasets of economic activity across Greater London. The \texttt{GLA\_data} folder 
contains the data compiled from the two GLA studies on ward profiles and 
retail floor space. The \texttt{dest\_sizes.csv} and \texttt{origin\_sizes.csv} are the destination and 
origin zone sizes used in the paper. The \texttt{exp\_times.csv} and \texttt{exp\_distances.csv} 
are the two different transport network metrics used, calculated via 
$\exp(-d_{ij}/\max(d_{ij}))$ from the respective \texttt{distances.csv} and \texttt{times.csv} files. The \texttt{Google\_Distance\_Matrix\_Data} folder contains transport times and distances using 
the Google Maps API service. Each file is a pickle-dictionary containing the API output for different 
travel modes: \texttt{transit} (public transport) and \texttt{driving} (driving, no traffic).
The \texttt{departure\_time} for transit is Sunday, June 19th 2022, 1 pm GMT (in Unix time: \texttt{departure\_time = 1655640000}). However, since
trips in the past cannot be computed, a future date must always be specified when using the 
API. The data is also available as a cost matrix in \texttt{.csv} format: entries are given in seconds and metres
respectively.
\subsection*{Running the code}
As with the SIR model, all the data (including the data from the GLA and Google Maps Distance API) is stored in the \texttt{data/HarrisWilson} folder. All data is given in \texttt{.csv} format. You can load data using the \texttt{load\_from\_dir} key:
\begin{lstlisting}
HarrisWilson:
  Data:
    load_from_dir:
      network: data/London_data/exp_times.csv
      origin_zones: data/London_data/origin_sizes.csv
      destination_zones: data/London_data/dest_sizes.csv
\end{lstlisting}
Here, we are specifying the exact locations of the network, origin zone and destination zone size files. You can also pass a single directory to \texttt{load\_from\_dir}, as long as that directory contains \texttt{.csv} files labelled as \texttt{origin\_sizes.csv}, \texttt{dest\_sizes.csv}, and \texttt{network.csv}:
\begin{lstlisting}
HarrisWilson:
  Data:
    load_from_dir: path/to/data
\end{lstlisting}    
Run the model using the corresponding configuration sets to reproduce the plots; for example, to produce the plots of the London dataset, just do
\begin{lstlisting}
utopya run HarrisWilson --cfg-set London_dataset
\end{lstlisting}

\end{document}